\documentclass[11pt]{amsart}
\usepackage[top=1.5in, bottom=1.5in, left=1.5in, right=1.5in]{geometry}
\expandafter\let\csname ver@amsthm.sty\endcsname\relax

\usepackage[all]{xy}
\usepackage{graphicx} 
\usepackage{tikz}
\usetikzlibrary{patterns}
\usetikzlibrary{decorations.pathreplacing}
\usetikzlibrary{arrows}
\usetikzlibrary{decorations.markings}

\usepackage{amsmath}
\usepackage{amssymb}
\usepackage{enumerate}
\usepackage{mathdots}
\usepackage{mathtools}

\usepackage{hyperref}
\usepackage{amsthm}
\usepackage[capitalize,noabbrev]{cleveref}

\allowdisplaybreaks

\usepackage{ytableau}
\ytableausetup{boxsize=1em}

\numberwithin{equation}{section}

\newtheorem{thm}{Theorem}[section]
\newtheorem{lemma}[thm]{Lemma}

\newtheorem{Definition}[thm]{Definition}

\newtheorem{Example}[thm]{Example}

\newtheorem{Remark}[thm]{Remark}
\newenvironment{remark}
  {\begin{Remark}\rm}{\end{Remark}}
  
\crefname{thm}{Theorem}{Theorems}
\crefname{lemma}{Lemma}{Lemmas}
\crefname{cor}{Corollary}{Corollaries}
\crefname{prop}{Proposition}{Propositions}
\crefname{conj}{Conjecture}{Conjectures}
\crefname{question}{Question}{Questions}
\crefname{problem}{Problem}{Problems}

\crefname{definition}{Definition}{Definitions}
\crefname{example}{Example}{Examples}
\crefname{remark}{Remark}{Remarks}
\newcommand{\M}{\operatorname{M}}

\newcommand{\emailhref}[1]{\email{\href{#1}{#1}}}

\newcommand{\dfn}[1]{\textcolor{blue}{\emph{#1}}}

\title[Plane partitions of shifted double staircase shape]{Plane partitions of \\ shifted double staircase shape}
\author[S. Hopkins]{Sam Hopkins}\emailhref{shopkins@umn.edu}
\address{School of Mathematics, University of Minnesota, Minneapolis, MN 55455}
\author[T. Lai]{Tri Lai}\emailhref{tlai3@unl.edu}
\address{Department of Mathematics, University of Nebraska, Lincoln, NE 68588}


\date{\today}

\begin{document}

\begin{abstract}
We give a product formula for the number of shifted plane partitions of shifted double staircase shape with bounded entries. This is the first new example of a family of shapes with a plane partition product formula in many years. The proof is based on the theory of lozenge tilings; specifically, we apply the ``free boundary'' Kuo condensation due to Ciucu.
\end{abstract}

\maketitle

\section{Introduction and statement of results} \label{sec:intro} 

An \dfn{$a \times b$ plane partition} is an $a\times b$ array $\pi=(\pi_{i,j})$ of nonnegative integers that is weakly decreasing along rows and down columns. Let $\mathcal{PP}^m(a\times b)$ denote the set of such plane partitions with largest entry less than or equal to $m$. MacMahon's celebrated product formula~\cite{macmahon1915combinatory} for the number of these plane partitions is
\[\#\mathcal{PP}^m(a\times b) = \prod_{i=1}^{a} \prod_{j=1}^{b} \frac{m+i+j-1}{i+j-1}.\]
MacMahon's formula is nowadays recognized as one of the most elegant in algebraic and enumerative combinatorics.

In this paper we consider plane partitions of other shapes beyond rectangles. For~$\lambda$ a partition, a \dfn{plane partition of (unshifted) shape $\lambda$} is a filling of the Young diagram of~$\lambda$ with nonnegative integers that is weakly decreasing along rows and down columns. Let $\mathcal{PP}^m(\lambda)$ denote the set of such plane partitions with largest entry less than or equal to $m$. Similarly, for $\lambda$ a strict partition, a \dfn{(shifted) plane partition of shifted shape $\lambda$} is a filling of the shifted Young diagram of $\lambda$ with nonnegative integers that is weakly decreasing along rows and down columns. Let $\mathcal{SPP}^m(\lambda)$ denote the set of such plane partitions with largest entry less than or equal to $m$.

Plane partitions of other unshifted/shifted shapes beyond rectangles have received a significant amount of attention, especially since the $1970$s, because of their connection to \dfn{symmetry classes} of plane partitions. For example, the \dfn{symmetric} plane partitions in $\mathcal{PP}^{m}(n \times n)$ are evidently in bijection with the plane partitions in $\mathcal{SPP}^{m}(\lambda)$ when $\lambda = (n,n-1,\ldots,1)$ is a \dfn{shifted staircase} shape. MacMahon conjectured~\cite{macmahon1899partitions}, and Andrews~\cite{andrews1977macmahon} proved, an elegant product formula for the number of symmetric plane partitions. We will review all known unshifted/shifted shapes with plane partition product formulas in \cref{sec:background} below.

Here we establish a plane partition product formula for a new family of shapes:

\begin{thm} \label{thm:sds}
Let $0\leq k \leq n$, and let $\lambda \coloneqq  (n,n-1,\ldots,1) + (k,k-1,\ldots,1)$ be a \dfn{shifted double staircase} shape. Then
\[\#\mathcal{SPP}^{m}(\lambda) = \prod_{1\leq i \leq j \leq n} \frac{m+i+j-1}{i+j-1} \prod_{1\leq i \leq j \leq k} \frac{m+i+j}{i+j}.\]
\end{thm}

\begin{figure}
\begin{center}
\begin{ytableau}
3 & 3 & 3 & 2 & 2 & 2 & 2 &2 &2\\
\none & 2 & 2 & 2 & 2 & 2 & 1 &1 &\none \\
\none & \none & 2 & 2 & 1 & 1 & 0 &\none &\none \\
\none & \none & \none & 1 & 1 & 0 & \none &\none &\none \\
\none & \none & \none &\none & 1 & 0 & \none &\none &\none \\
\none & \none & \none & \none & \none & 0 & \none &\none &\none
\end{ytableau}
\end{center}
\caption{A plane partition of shifted double staircase shape with $k=3$, $n=6$, $m=3$.} \label{fig:sds_ex}
\end{figure}

\Cref{fig:sds_ex} depicts a plane partition of shifted double staircase shape in the case $k=3$, $n=6$, $m=3$.

The case $k=0$ of \cref{thm:sds} corresponds to shifted staircases, where the formula is known by the aforementioned result of Andrews concerning symmetric plane partitions. The cases $k=n-1$ and $k=n$ of \cref{thm:sds} are also known: these shapes are \dfn{shifted trapezoids}, and Proctor~\cite{proctor1983trapezoid} gave a product formula for their plane partition enumeration. In fact, the cases $k=n-1$ and $k=n$ correspond to a symmetry class of plane partitions: the \dfn{symmetric, self-complementary} plane partitions. The other cases of \cref{thm:sds} beyond $k=0$, $k=n-1$, and $k=n$ are new. This is the first new example of a family of shapes with a plane partition product formula in many years.

\begin{figure}
\begin{center}
\begin{tabular}{c c c}
\parbox{7.5cm}{\includegraphics[width=7.5cm]{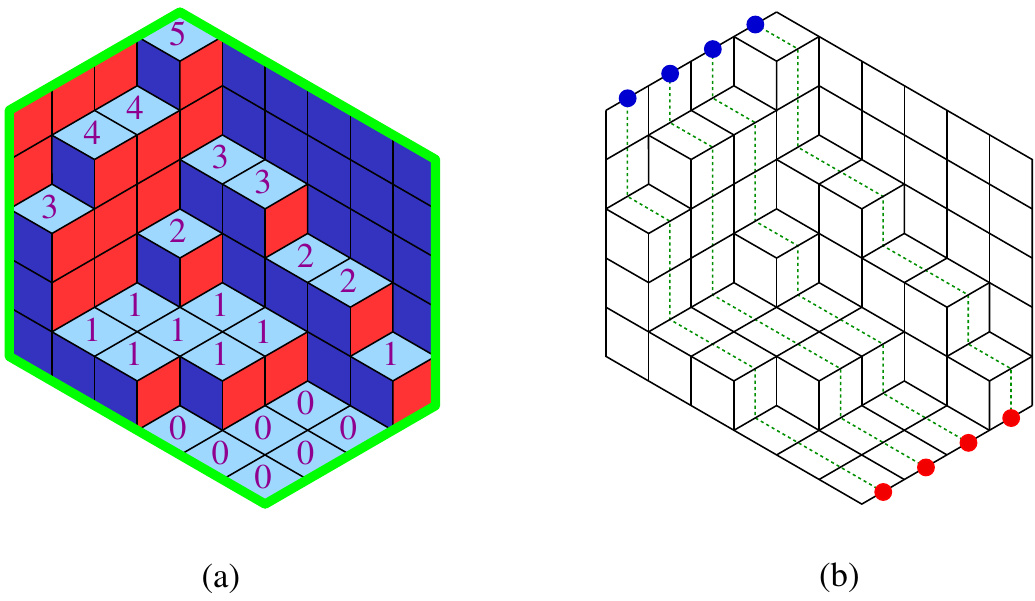}}  & \Large $\Leftrightarrow$ 
 & \parbox{1in}{\begin{ytableau}
5 & 3 & 3 & 2 & 2 & 1\\
4 & 2 & 1 & 1 & 1 & 0 \\
4 & 1 & 1 & 1 & 0 & 0 \\
3 & 1 & 1 & 0 & 0 & 0
\end{ytableau}}  \end{tabular}
\end{center}
\caption{The correspondence between boxed plane partitions and lozenge tilings of a hexagon. Each row of the plane partition records the label of the horizontal lozenges in picture~(a) along a lattice path in picture ~(b).} \label{fig:pparts_to_lozenges}
\end{figure}

Our proof of \cref{thm:sds} is based on the theory of \dfn{lozenge tilings of the triangular lattice}. There is a well-known bijective correspondence between the plane partitions in~$\mathcal{PP}^{m}(a\times b)$ and the lozenge tilings of a hexagonal region (with sides of lengths $a$, $b$, and $m$) of the triangular lattice. Briefly, we draw $\pi \in \mathcal{PP}^{m}(a\times b)$ as a stacking of cubes in an $a\times b \times m$ box, with $\pi_{i,j}$ cubes stacked at position $(i,j)$; then we ``forget'' that the picture we drew represented a three-dimensional object, and instead view it as a flat, two-dimensional object. This is depicted in \cref{fig:pparts_to_lozenges}. 

\begin{figure}\centering
\includegraphics[width=3in]{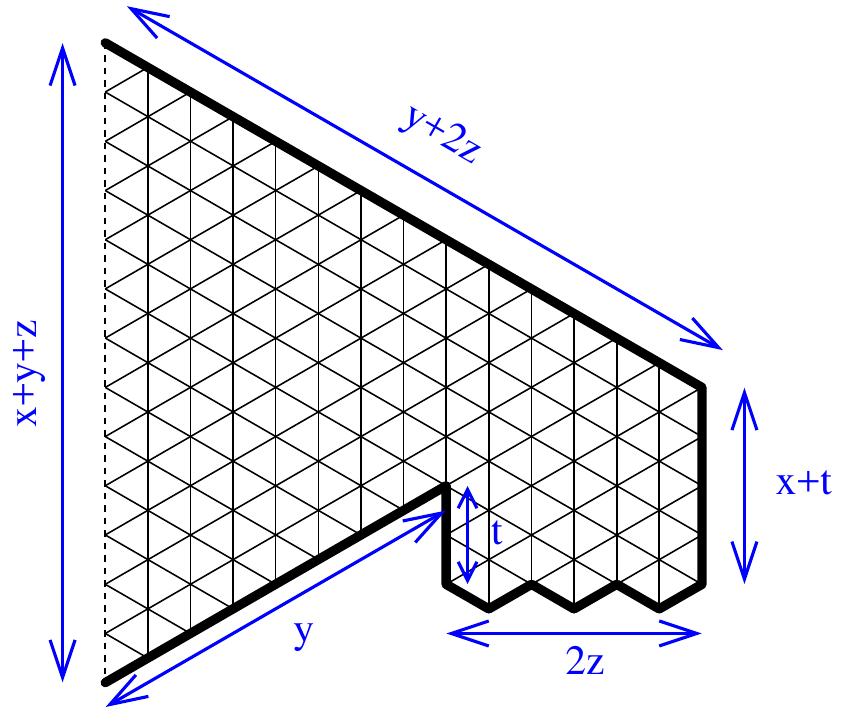}
\caption{The flashlight region $F_{x,y,z,t}$. The dotted side indicates the free boundary.} \label{fig:flashlight}
\end{figure}

Plane partitions of other shapes similarly correspond to lozenge tilings of other kinds of regions. For our purposes, the relevant region is the \dfn{flashlight region} $F_{x,y,z,t}$ depicted in \cref{fig:flashlight}. In this region, the dashed line represents a \dfn{free boundary} from which lozenges are allowed to protrude. Free boundaries appear in regions which correspond to plane partitions of shifted shapes (or equivalently, symmetric plane partitions of unshifted shapes).

\begin{figure}
\begin{center}
\begin{tabular}{c c c}
  \parbox{7cm}{\includegraphics[width=7cm]{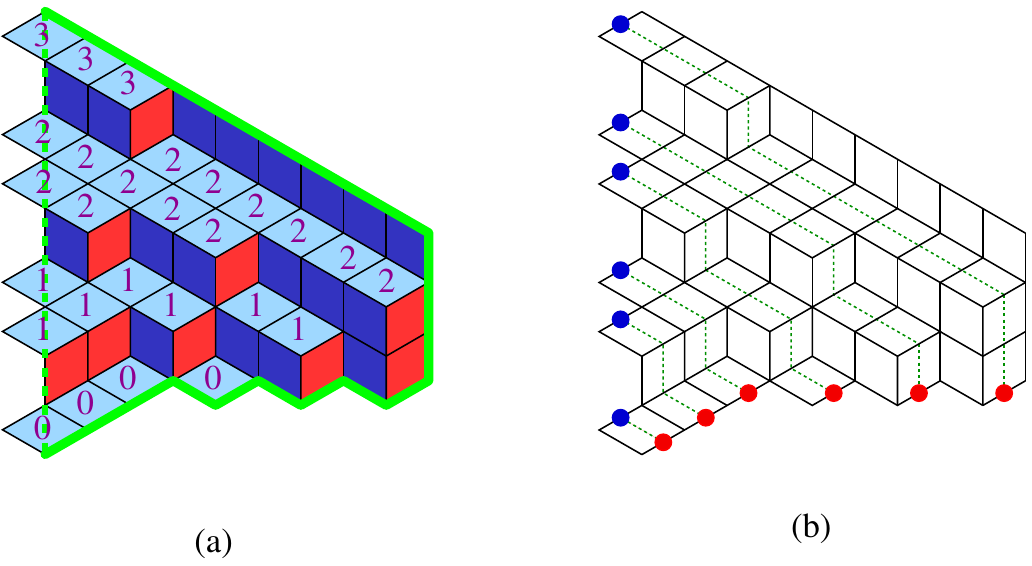}}  & \Large $\Leftrightarrow$ &\parbox{1in}{\begin{ytableau}
3 & 3 & 3 & 2 & 2 & 2 & 2 &2 &2\\
\none & 2 & 2 & 2 & 2 & 2 & 1 &1 &\none \\
\none & \none & 2 & 2 & 1 & 1 & 0 &\none &\none \\
\none & \none & \none & 1 & 1 & 0 & \none &\none &\none \\
\none & \none & \none &\none & 1 & 0 & \none &\none &\none \\
\none & \none & \none & \none & \none & 0 & \none &\none &\none
\end{ytableau}}
\end{tabular}
\end{center}
\caption{The correspondence between the lozenge tilings of the flashlight region (with $t=0$) and plane partitions of shifted double staircase shape. Each row of the plane partition records the label of the horizontal lozenges in picture~(a) along a lattice path in picture ~(b).} \label{fig:correspondence}
\end{figure}

We prove the following product formula for the enumeration of lozenge tilings of the regions $F_{x,y,z,t}$:

\begin{thm} \label{thm:flashlight}
For nonnegative integers $x$, $z$, $t$ and positive integer $y$, the number of lozenge tilings of $F_{x,y,z,t}$ is given by
\begin{align}\label{mainformula}
\M(F_{x,y,z,t})= \hspace{-10pt} \prod_{1\leq i \leq j \leq y+z} \hspace{-10pt} \frac{x+i+j-1}{i+j-1} \; \; \prod_{1\leq i \leq j \leq z} \hspace{-5pt} \frac{x+i+j}{i+j} \; \; \; \prod_{i=1}^{t} \prod_{j=1}^{z} \frac{x+z+2i+j}{x+2i+j-1},
\end{align}
where we use the notation $\M(R)$ for the number of lozenge tilings of the region $R$ and where empty products are taken to be $1$ by convention.
\end{thm}

The case $t=0$ of \cref{thm:flashlight} yields \cref{thm:sds}, where we take $x=m$, $y+z=n$, and $z=k$: see \cref{fig:correspondence}.

Lozenge tilings of the flashlight region were previously considered, and enumerated, by Ciucu in~\cite{ciucu2020correlation}. However, the regions addressed there are a subset of the ones we consider here: e.g., Ciucu's regions essentially had the requirement $t\geq \lfloor y/2 \rfloor$. (We remark that our parameterization of these regions is slightly different than Ciucu's, and it leads to a ``cleaner,'' although equivalent, product formula even in the cases he considered.) The case $t=0$, the one relevant for the shifted double staircase plane partitions, therefore was not directly addressed in~\cite{ciucu2020correlation}. Nevertheless, the techniques that Ciucu developed in~\cite{ciucu2020correlation} are sufficient for proving \cref{thm:flashlight}.

So what are these techniques? In 2003, Kuo~\cite{kuo2004applications} invented a powerful method for proving product formulas for lozenge tilings and other \dfn{dimer} (or \dfn{perfect matching}) enumeration problems, which is now called \dfn{Kuo condensation}. Kuo condensation is essentially a recurrence relation satisfied by these dimer counts, and often allows one to easily prove product formulas, provided one can ``guess'' the product formula in advance. In~\cite{ciucu2020correlation}, Ciucu established a variant of Kuo condensation that works for lozenge tilings of regions with a free boundary. A straightforward application of this free boundary Kuo condensation yields \cref{thm:flashlight}.

Before we finish this introduction, let us briefly review the history behind the discovery of \cref{thm:sds}. Plane partitions for a given unshifted/shifted shape $\lambda$ are a special case of the more general notion of \dfn{$P$-partitions} for partially ordered sets~$P$. In particular, the enumerations $\#\mathcal{PP}^{m}(\lambda)$ and $\#\mathcal{SPP}^{m}(\lambda)$ are special cases of \dfn{order polynomials} of posets. In~\cite{hopkins2020order}, the first author proposed an apparently powerful heuristic which says that the posets which have product formulas for their order polynomials are the same as the posets with ``good dynamical behavior.'' Here ``dynamical behavior'' means the behavior of various invertible operators coming from algebra on objects associated to the poset, such as Sch\"{u}tzenberger's \dfn{promotion of linear extensions}. In~\cite{haiman1992dual, haiman1992characterization}, Haiman showed that the shapes with good behavior of promotion are exactly: rectangles, staircases, shifted trapezoids, and shifted double staircases. It was already known (see \cref{sec:background}) that rectangles, staircases, and shifted trapezoids have product formulas for their order polynomials, i.e., product formulas enumerating their plane partitions. \Cref{thm:sds} shows the same is true for the shifted double staircases. Indeed, it is precisely via this heuristic that \cref{thm:sds} was discovered.

The rest of the paper is structured as follows: in \cref{sec:background} we go over more background on counting plane partitions of different shapes; and in \cref{sec:proof} we prove \cref{thm:flashlight} (and thus \cref{thm:sds}) via the free boundary Kuo condensation.

\begin{remark} \label{rem:okada}
After the first version of this article was available online, Okada~\cite{okada2020intermediate} gave another proof of~\cref{thm:sds}. Actually, he gave an algebraic extension of~\cref{thm:sds} involving characters of classical groups and their variants, especially the \dfn{intermediate symplectic group characters} of Proctor~\cite{proctor1988odd}. In fact, Okada proved an algebraic extension of the four parameter product formula in~\cref{thm:flashlight} for the number of tilings of the flashlight region. Okada's algebraic extensions in particular yield $q$-analogs of these product formulas.
\end{remark}

\medskip

\noindent {\bf Acknowledgments}: S.H. thanks Robert Proctor for useful comments, and thanks Soichi Okada for sharing his work~\cite{okada2020intermediate} while it was still in progress and conjectural. We also thank the anonymous referees for their careful attention to our paper and useful comments. S.H. was supported by NSF grant $\#1802920$. T.L. was supported by Simons Foundation Collaboration Grant $\#585923$.

\section{Background on counting plane partitions of different shapes} \label{sec:background}

In this section we review the known results concerning plane partition enumeration formulas for different shapes. For a more complete account of the history of plane partitions, with an emphasis on symmetry classes, see~\cite{krattenthaler2016plane}. Nothing in this section is necessary for the proof of the main result, so the reader who is not interested in or already knows the context may safely skip it.

Let us first formalize notions introduced in \cref{sec:intro}. For $\lambda=(\lambda_1,\lambda_2,\ldots)$ a \dfn{partition}, its \dfn{Young diagram} is the left-justified collection of boxes which has $\lambda_i$ boxes in row $i$.  As we said earlier, a \dfn{plane partition of (unshifted) shape $\lambda$} is a filling of the Young diagram of~$\lambda$ with nonnegative integers that is weakly decreasing along rows and down columns; and we use $\mathcal{PP}^m(\lambda)$ to denote the set of such plane partitions with largest entry less than or equal to $m$. 

\begin{figure}
\begin{center}
\begin{tabular}{c c c}
\parbox{1in}{\begin{ytableau}
5 & 3 & 3 & 2 & 2 \\
4 & 2 & 1 & 1  \\
4 & 1 & 1 & 1  \\
3 & 1 & 0
\end{ytableau}} & \Large $\Leftrightarrow$ & \parbox{1.5in}{\includegraphics[width=1.25in]{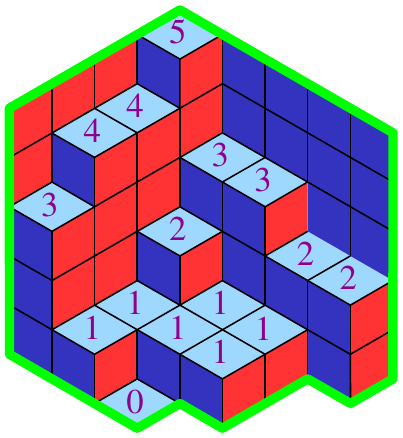}}
\end{tabular}
\end{center}
\caption{The correspondence between plane partitions and lozenge tilings for arbitrary unshifted shapes.} \label{fig:nonrect_pparts_to_lozenges}
\end{figure}

A \dfn{lozenge} is a rhombus consisting of two adjacent (i.e., sharing a common edge) unit triangles of the triangular lattice. There are three possible ``orientations'' of lozenges, which we often depict as three different colors. A \dfn{lozenge tiling} of a region of the triangular lattice is a covering of this region with non-overlapping lozenges. As we explained above (see \cref{fig:pparts_to_lozenges}) there is a bijective correspondence between $\mathcal{PP}^{m}(a \times b)$ and lozenge tilings of a hexagonal region of the triangular lattice: we first represent the plane partition as 3D stack of cubes; then we treat the 3D picture as a 2D picture. (While the 3D picture of plane partitions goes all the way back to MacMahon~\cite{macmahon1899partitions}, the connection to lozenge tilings is much more recent~\cite{guy1989calissons}.) For other unshifted shapes $\lambda$ beyond rectangles, there is again a bijective correspondence between $\mathcal{PP}^{m}(\lambda)$ and lozenge tilings of a region of the triangular lattice, where the region is obtained from a hexagon by appropriately modifying its ``lower boundary'' to follow the shape of $\lambda$: see \cref{fig:nonrect_pparts_to_lozenges}.

For $\lambda$ a \dfn{strict partition}, its \dfn{shifted Young diagram} is almost the same as its Young diagram, except that we indent each row one position to the right compared to the row above it. As we said earlier, a \dfn{(shifted) plane partition of shifted shape $\lambda$} is a filling of the shifted Young diagram of~$\lambda$ with nonnegative integers that is weakly decreasing along rows and down columns; and we use $\mathcal{SPP}^m(\lambda)$ to denote the set of such plane partitions with largest entry less than or equal to $m$.

For $\lambda$ a strict partition, we use $\lambda^D$ to denote the partition whose Young diagram is obtained from the shifted Young diagram of $\lambda$ by gluing a reflected copy of that shifted Young diagram along the main diagonal. Note that such a $\lambda^D$ is always \dfn{self-conjugate} (i.e., invariant under reflection across the main diagonal). If $\mu$ is any self-conjugate partition, then we use $\mathrm{Tr}\colon\mathcal{PP}^{m}(\mu)\to\mathcal{PP}^{m}(\mu)$ to denote the \dfn{transposition} map (i.e., reflection across the main diagonal) acting on plane partitions of shape~$\mu$. There is an obvious bijection between $\mathcal{SPP}^m(\lambda)$ and the set of \dfn{symmetric} (i.e., transposition invariant) plane partitions in~$\mathcal{PP}^{m}(\lambda^D)$.

\begin{figure}
\begin{center}
\begin{tabular}{c c c c c c c}
\parbox{0.75in}{\begin{ytableau}
4 & 4 & 3 & 2 & 1 \\
\none & 3 & 3 & 2 & 0  \\
\none & \none & 2 & 1  \\
\none & \none & \none & 1
\end{ytableau}} & \Large $\Leftrightarrow$ &  \parbox{0.8in}{\begin{ytableau}
4 & 4 & 3 & 2 & 1 \\
4 & 3 & 3 & 2 & 0  \\
3 & 3 & 2 & 1  \\
2 & 2 & 1 & 1 \\
1 & 0
\end{ytableau}} & \Large $\Leftrightarrow$ & \parbox{1.25in}{\includegraphics[height=1.25in]{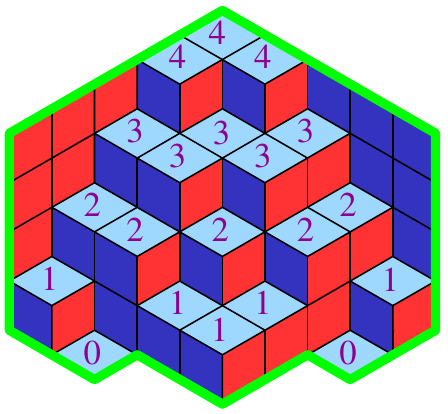}} & \Large $\Leftrightarrow$ & \parbox{0.75in}{\includegraphics[height=1.25in]{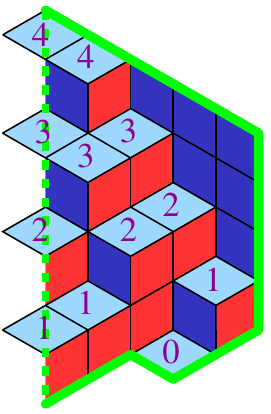}}
\end{tabular}
\end{center}
\caption{The correspondence between shifted plane partitions and lozenge tilings of regions with a free boundary, going through symmetric plane partitions/lozenge tilings. In the right picture the free boundary is dashed.} \label{fig:shifted_pparts_to_lozenges}
\end{figure}

The above discussion tells us that plane partitions in $\mathcal{SPP}^m(\lambda)$ are in bijection with lozenge tilings of the region of the triangular lattice corresponding to $\lambda^D$ which are \dfn{symmetric} (i.e., invariant under reflection across the vertical axis of symmetry). By only remembering the ``right half'' of a such a symmetric lozenge tiling, we obtain a lozenge tiling of a ``halved'' region, where the symmetry axis becomes a \dfn{free boundary}. In a lozenge tiling of a region with a free boundary, lozenges are allowed to protrude from the free boundary. The upshot is that plane partitions in $\mathcal{SPP}^m(\lambda)$ are in bijection with lozenge tilings of a region of the triangular lattice with a free boundary. This whole construction is depicted in \cref{fig:shifted_pparts_to_lozenges}.

We have now seen that the plane partitions of every unshifted shape (respectively, shifted shape) correspond to lozenge tilings of a region of the triangular lattice (resp., region with free boundary). Importantly, the converse is not true; that is, there are many regions whose tilings do not correspond to plane partitions. Indeed, this is part of what makes the lozenge tilings perspective useful: regions which do not correspond to plane partitions naturally arise in recurrences even when one is only interested in the plane partition case.

Now let us briefly review the main techniques for counting plane partitions. There are three such main techniques:
\begin{itemize}
\item {\bf determinantal or Pfaffian formulas} coming from {\bf nonintersecting lattice paths} interpretations;
\item formulas deduced from the {\bf representation theory of classical groups and their Lie algebras};
\item {\bf recurrence relations} satisfied by these counts, like {\bf Kuo condensation} and its variants.
\end{itemize}
Of course, these three techniques have a lot of overlap. For instance, the \dfn{Weyl character formula} expresses the character of an irreducible representation of a classical group as a determinant. And Kuo condensation is really a combinatorialization of \dfn{Dodgson condensation}/the \dfn{Desnanot--Jacobi identity} for computing determinants. Furthermore, there are other slightly different techniques which have been successfully used to count plane partitions: e.g., Kuperberg~\cite{kuperberg1998permanent} championed the use of Kasteleyn's \dfn{permanent-determinant}/\dfn{Hafnian-Pfaffian} method from dimer theory. Kasteleyn's method is similar to, but slightly different from, the nonintersecting lattice paths technique.

In the remainder of this section we will see how these techniques apply to yield product formulas in different special cases. But first let us discuss formulas for arbitrary shapes.

\subsection{Formulas for arbitrary shapes}

For unshifted shapes $\lambda$ we have the following \dfn{determinantal} theorem of MacMahon:

\begin{thm}[{MacMahon~\cite{macmahon1915combinatory}}] \label{thm:macmahon_determinant}
Let $\lambda=(\lambda_1,\lambda_2,\ldots,\lambda_n)$ be a partition. Then
\[ \#\mathcal{PP}^{m}(\lambda) = \mathrm{det} \left ( \binom{\lambda_i+m}{i-j+m} \right )_{i,j=1,\ldots,n}.\]
Here $\binom{a}{b}=0$ if $b < 0$.
\end{thm}

Actually, MacMahon~\cite[Volume 2, \S X, Chapter II]{macmahon1915combinatory} proved a \dfn{$q$-analog} of \cref{thm:macmahon_determinant}. He showed that
\[ \sum_{\pi \in \mathcal{PP}^{m}(\lambda)}q^{|\pi|} =  \mathrm{det} \left ( q^{f(i,j)} \binom{\lambda_i+m}{i-j+m}_q \right )_{i,j=1,\ldots,n}, \]
where $|\pi| \coloneqq  \sum_{u \in \lambda} \pi(u)$ is the \dfn{size} of $\pi$, i.e., the sum of the entries of the plane partition, and $f(i,j) \coloneqq  \begin{cases} \binom{j-i}{2} &\textrm{if $j>i$}; \\ \binom{i-j+1}{2} &\textrm{if $i \geq j$}\end{cases}$, and $\binom{a}{b}_q \coloneqq  \prod_{1\leq i \leq b} \frac{(1-q^{a+1-i})}{(1-q^i)}$ is the usual \dfn{$q$-binomial coefficient}. It is also worth noting that Kreweras~\cite{kreweras1965classe} proved a generalization of \cref{thm:macmahon_determinant} which works for \dfn{skew shapes} $\lambda/\mu$.

A modern approach to \cref{thm:macmahon_determinant} is to encode the plane partitions in $\mathcal{PP}^{m}(\lambda)$ as collections of \dfn{nonintersecting lattice paths} connecting a particular set of sources to a particular set of sinks on the rectangular grid, and then to apply the \dfn{Lindstr\"{o}m--Gessel--Viennot} lemma to count these nonintersecting lattice paths. See the seminal preprint~\cite{gessel1988determinants} for the details. 

The natural encoding of plane partitions of shifted shapes via nonintersecting lattice paths involves paths whose set of sources are fixed, but whose set of sinks are allowed to vary in some set. (This is basically the same phenomenon as the free boundary in the lozenge tilings corresponding to such plane partitions and is visible for instance in \cref{fig:correspondence}.) Stembridge~\cite{stembridge1990pfaffian}, building on a theorem of Okada~\cite{okada1989generating} about the sum of all minors of a matrix, developed the relevant machinery for counting these kind of nonintersecting lattice paths with variable endpoints: rather than a determinant, he showed that the answer can be expressed as a \dfn{Pfaffian} $\mathrm{Pf}(M)$ of a matrix $M$. In the case of plane partitions of shifted shapes, this method gives:

\begin{thm}[{Stembridge~\cite{stembridge1990pfaffian}}] \label{thm:stembridge_pfaffian}
Let $\lambda=(\lambda_1,\lambda_2,\ldots,\lambda_n)$ be a strict partition, and assume $n$ is even by setting $\lambda_n \coloneqq  0$ if necessary. Let $m\in \mathbb{N}$. Let $M$ be the $n\times n$ skew-symmetric matrix whose $(i,j)$ entry for $i,j=1,\ldots,n$ is:
\begin{align*}
M_{i,j}=\sum_{0\leq k \leq \ell \leq m+n-1} &\left [ \binom{\lambda_i-1+m+i-1-k}{m+i-1-k}\binom{\lambda_j-1+m+j-1-\ell}{m+j-1-\ell} \right. \\
&\left. - \binom{\lambda_i-1+m+i-1-\ell}{m+i-1-\ell}  \binom{\lambda_j-1+m+j-1-k}{m+j-1-k} \right ].
\end{align*}
Here $\binom{a}{b}=0$ if $b < 0$ (but $\binom{-1}{0}=1$). Then $\# \mathcal{SPP}^{m}(\lambda) = \mathrm{Pf} (M)$.
\end{thm}

Note that Stembridge never explicitly stated \cref{thm:stembridge_pfaffian}, but his techniques easily allow one to prove it (and see~\cite[Proof of Theorem~5.1]{stembridge1990pfaffian} for the basic idea behind the encoding). In fact, it is again possible to upgrade  \cref{thm:stembridge_pfaffian} to a $q$-analog which gives the size generating function for these shifted plane partitions.

As one can see, the entries of the matrix in \cref{thm:stembridge_pfaffian} are much more complicated than the entries of the matrix in \cref{thm:macmahon_determinant}, and this reflects the general principle that shifted shapes are harder to deal with than unshifted shapes. As far as we are aware, \cref{thm:macmahon_determinant,thm:stembridge_pfaffian} are the best general formulas for counting $\#\mathcal{PP}^{m}(\lambda)$ or $\#\mathcal{SPP}^{m}(\lambda)$ for arbitrary $\lambda$.

\subsection{Product formulas for specific shapes}

\begin{table}
\ytableausetup{boxsize=0.75em}
\begin{center}
\renewcommand{\arraystretch}{1.5}
\noindent
\makebox[\textwidth]{{
\begin{tabular}{c | c | c | c | c }
Shape & Product formula & \parbox{1.2in}{\begin{center} Lozenge tiling region\end{center}} & \parbox{1.2in}{\begin{center}Plane partition symmetry class\end{center}} & References \\ \hline
\parbox{1.75in} {\begin{center} Rectangle \\
\begin{tikzpicture}
\node at (0,0){\ydiagram{10,10,10,10}};
\draw [decorate,decoration={brace,amplitude=5pt},yshift=0pt] (-1.5,0.7) -- (1.5,0.7) node [black,midway,yshift=0.4cm]  {$b$};
\draw [decorate,decoration={brace,amplitude=5pt,mirror},yshift=0pt] (-1.6,0.6) -- (-1.6,-0.6) node [black,midway,xshift=-0.4cm]  {$a$};
\end{tikzpicture}\end{center}} & $ \displaystyle \prod_{i=1}^{a} \prod_{j=1}^{b} \frac{m+i+j-1}{i+j-1}$ & \parbox{1in}{\includegraphics[width=1.1in]{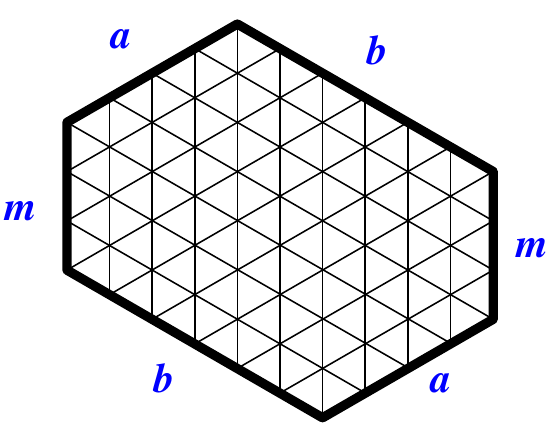}} & All & \parbox{1in}{\begin{center} MacMahon~\cite{macmahon1915combinatory}; {\bf Extensions}: holes~\cite{ciucu1998enumeration, helfhott1999enumeration, okada1998number, ciucu2001enumeration, ciucu2013dual, ciucu2017other, lai2020lozenge}; dents~\cite{eisenkolbl1999rhombus, ciucu2015generalization, lai2017qenumeration} \end{center}}\\ \hline
\parbox{1.75in} {\begin{center} Shifted staircase \\
\begin{tikzpicture}
\node at (0,0){\ydiagram{7,1+6,2+5,3+4,4+3,5+2,6+1}};
\draw [decorate,decoration={brace,amplitude=5pt},yshift=0pt] (-1.1,1.1) -- (1.1,1.1) node [black,midway,yshift=0.4cm]  {$n$};
\end{tikzpicture}\end{center}} & $ \displaystyle \prod_{1\leq i \leq j \leq n} \frac{m+i+j-1}{i+j-1}$ & \parbox{1in}{\includegraphics[width=0.8in]{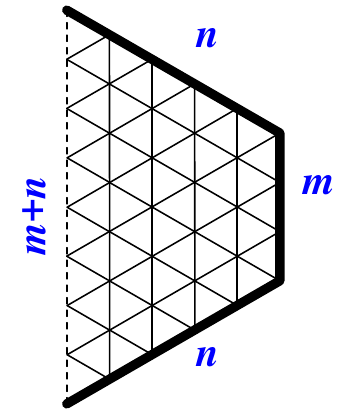}} & Symmetric & \parbox{1.2in}{\vspace{.5\baselineskip} \begin{center} Andrews~\cite{andrews1977macmahon, andrews1978plane, andrews1977plane}, Macdonald~\cite{macdonald1979symmetric}, Gordon~\cite{gordon1983bender}, Proctor~\cite{proctor1984bruhat, proctor1990new}; {\bf Extensions}: see shifted double staircase below \end{center} \vspace{.5\baselineskip}} \\ \hline
\parbox{1.75in} {\begin{center} Staircase \\
\begin{tikzpicture}
\node at (0,0){\ydiagram{10,9,8,7}};
\draw [decorate,decoration={brace,amplitude=5pt},yshift=0pt] (-1.5,0.7) -- (1.5,0.7) node [black,midway,yshift=0.4cm]  {$b$};
\draw [decorate,decoration={brace,amplitude=5pt,mirror},yshift=0pt] (-1.6,0.6) -- (-1.6,-0.6) node [black,midway,xshift=-0.4cm]  {$a$};
\end{tikzpicture}\end{center}} & \parbox{1.6in}{\vspace{.5\baselineskip} \begin{center} $ \displaystyle \prod_{i=1}^{a} \left [ \prod_{j=1}^{b-a+1} \frac{m+i+j-1}{i+j-1} \cdot \right.$ \\ $ \displaystyle  \left. \prod_{j=b-a+2}^{b-a+i} \frac{2m+i+j-1}{i+j-1} \right ]$ \end{center} \vspace{.5\baselineskip}} & \parbox{1in}{\includegraphics[width=1.1in]{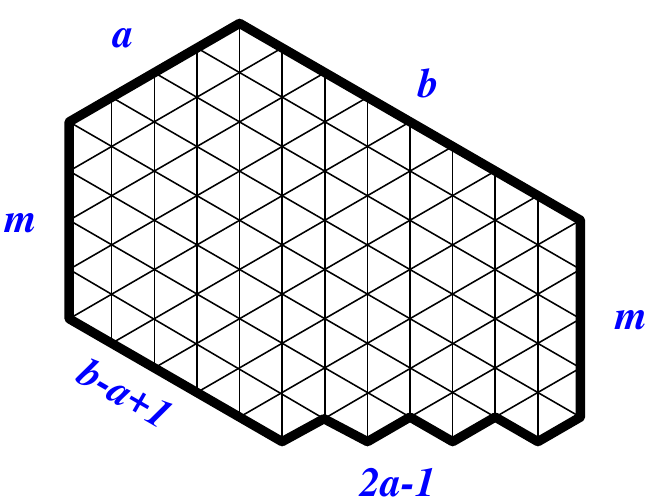}} & \parbox{1.2in}{\vspace{.5\baselineskip} \begin{center} Transpose-complementary $(a=b)$ \end{center} \vspace{.5\baselineskip}} & \parbox{1.3in}{\begin{center}Proctor~\cite{proctor1984research, proctor1988odd, proctor1990new}; {\bf Extensions}: Ciucu-Krattenthaler~\cite{ciucu2002enumeration} \end{center}}\ \\ \hline
\parbox{1.75in} {\begin{center} Shifted trapezoid \\
\begin{tikzpicture}
\node at (0,0){\ydiagram{10,1+8,2+6,3+4}};
\draw [decorate,decoration={brace,amplitude=5pt},yshift=0pt] (-1.5,0.7) -- (1.5,0.7) node [black,midway,yshift=0.4cm]  {$n$};
\draw [decorate,decoration={brace,amplitude=5pt,mirror},yshift=0pt] (-1.6,0.6) -- (-1.6,-0.6) node [black,midway,xshift=-0.4cm]  {$k$};
\end{tikzpicture}\end{center}} & $ \displaystyle \prod_{i=1}^{k}\prod_{j=1}^{n-k+1}  \frac{m+i+j-1}{i+j-1}$ & \parbox{1in}{\includegraphics[width=1.1in]{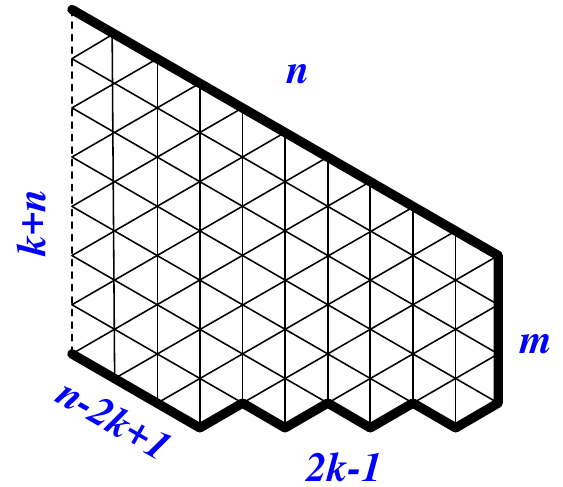}} & \parbox{1.2in}{\vspace{.5\baselineskip} \begin{center} Symmetric, self-complementary $(k=\lceil n/2 \rceil)$ \end{center} \vspace{.5\baselineskip}} & \parbox{1in}{\begin{center} Proctor~\cite{proctor1983trapezoid}; {\bf Extensions}: Ciucu~\cite{ciucu2015lozenge} \end{center}} \\ \hline
\parbox{1.75in} {\begin{center} Shifted double staircase \\
\begin{tikzpicture}
\node at (0,0){\ydiagram{10,1+8,2+6,3+4,4+3,5+2,6+1}};
\draw [decorate,decoration={brace,amplitude=5pt},yshift=0pt] (-1.5,1.1) -- (0.6,1.1) node [black,midway,yshift=0.4cm]  {$n$};
\draw [decorate,decoration={brace,amplitude=5pt},yshift=0pt] (0.6,1.1) -- (1.5,1.1) node [black,midway,yshift=0.4cm]  {$k$};
\end{tikzpicture} \end{center}} & \parbox{1.5in}{\begin{center}$ \displaystyle  \prod_{1\leq i \leq j \leq n} \frac{m+i+j-1}{i+j-1} \cdot $\\ $\displaystyle \prod_{1\leq i \leq j \leq k} \frac{m+i+j}{i+j}$\end{center}} & \parbox{1in}{\includegraphics[width=1.1in]{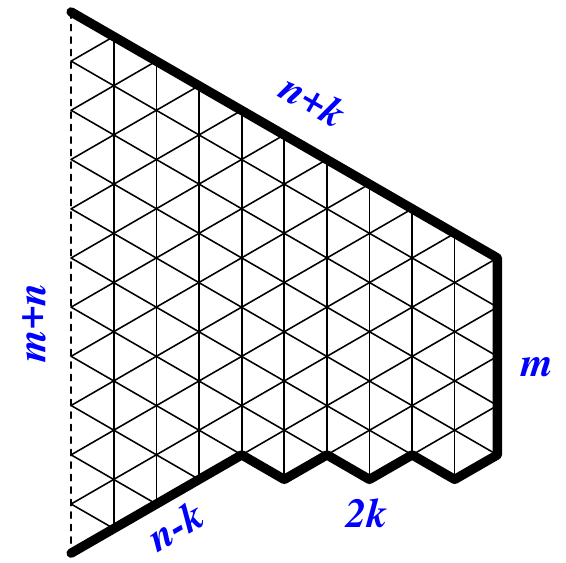}} & \parbox{1.2in}{\vspace{.5\baselineskip} \begin{center} Symmetric $(k=0)$; \\ Symmetric, self-complementary $(k=n-1,n)$ \end{center} \vspace{.5\baselineskip}} & \parbox{1.05in}{\begin{center} \Cref{thm:sds}; {\bf Extensions}: \cref{thm:flashlight} (c.f. Ciucu~\cite{ciucu2020correlation}) \end{center}}
\end{tabular}}}
\end{center}
\ytableausetup{boxsize=normal}
\caption{A table summarizing the known shapes with product formulas enumerating their plane partitions with entries $\leq m$. In the pictures of regions, the dashed lines represent free boundaries.} \label{table:shapes_table}
\end{table}

Now we review all the special families of unshifted/shifted shapes which have product formulas for their plane partition enumeration. These are all recorded in~\cref{table:shapes_table}. This table lists the product formula for the relevant $\#\mathcal{PP}^{m}(\lambda)$ or $\#\mathcal{SPP}^{m}(\lambda)$, and also shows the region of the triangular lattice whose lozenge tilings correspond to these plane partitions. It also shows the various \dfn{symmetry classes} of boxed plane partitions which are in bijection with plane partitions of these shapes. Let us take a moment to review these symmetry classes.

 Recall, as first outlined in~\cite{stanley1986symmetries}, that there are 10 classical symmetry classes of \dfn{boxed plane partitions}, i.e., plane partitions in $\mathcal{PP}^{m}(a\times b)$. (Many of these classes require $a=b$ or $a=b=m$ to exist.) Four of these symmetry classes are easily seen to be in bijection with plane partitions of certain shapes. To generate these four relevant symmetry classes, in addition to the already discussed transposition map, we also need the \dfn{complementation} map $\mathrm{Co}\colon \mathcal{PP}^{m}(a \times b) \to \mathcal{PP}^{m}(a \times b)$ which is defined by $\mathrm{Co}(\pi)_{i,j} \coloneqq  m - \pi_{a+1-i,b+1-j}$. (In the 3D picture of plane partitions, complementation is set-theoretic complementation inside of the $a\times b \times m$ box.) The four relevant symmetry classes of boxed plane partitions are then:
 \begin{itemize}
 \item \dfn{all} the plane partitions in $\mathcal{PP}^{m}(a \times b)$;
 \item the \dfn{symmetric} plane partitions in $\mathcal{PP}^{m}(n \times n)$, i.e., the plane partitions $\pi \in \mathcal{PP}^{m}(n \times n)$ with $\mathrm{Tr}(\pi)=\pi$;
 \item the \dfn{transpose-complementary} plane partitions in $\mathcal{PP}^{m}(n \times n)$, i.e., the plane partitions $\pi \in \mathcal{PP}^{m}(n \times n)$ with $\mathrm{Tr}(\pi)=\mathrm{Co}(\pi)$;
 \item the \dfn{symmetric, self-complementary} plane partitions in $\mathcal{PP}^{m}(n \times n)$, i.e., the plane partitions $\pi \in \mathcal{PP}^{m}(n \times n)$ with $\mathrm{Tr}(\pi)=\pi$ and $\mathrm{Co}(\pi)=\pi$.
 \end{itemize}  
 
 Let us now go through the rows of~\cref{table:shapes_table} in turn.

\subsubsection{Rectangles}

As mentioned at the beginning of \cref{sec:intro}, MacMahon~\cite{macmahon1915combinatory} famously proved that plane partitions of \dfn{rectangle shape} are counted by the product formula
\[\#\mathcal{PP}^m(a\times b) = \prod_{i=1}^{a} \prod_{j=1}^{b} \frac{m+i+j-1}{i+j-1}.\]
In fact, he established the following $q$-analog:
\[ \sum_{\pi \in \mathcal{PP}^m(a\times b)} q^{|\pi|} = \prod_{i=1}^{a} \prod_{j=1}^{b} \frac{(1-q^{m+i+j-1})}{(1-q^{i+j-1})}.\]

A modern representation theoretic explanation of MacMahon's product formula goes as follows. For a simple Lie algebra $\mathfrak{g}$ and an integral, dominant weight~$\lambda$ of~$\mathfrak{g}$, let~$V_{\mathfrak{g}}^{\lambda}$ denote the irreducible, finite-dimensional representation of $\mathfrak{g}$ with highest weight~$\lambda$. It is well-known (see, e.g.,~\cite{proctor1984bruhat, proctor1990new}) that the plane partitions in~$\mathcal{PP}^m(a\times b)$ index a basis of $V_{\mathfrak{gl}(a+b)}^{m^a}$. MacMahon's product formula for $\#\mathcal{PP}^m(a\times b)$ can then be deduced by applying the \dfn{Weyl dimension formula} to this representation $V_{\mathfrak{gl}(a+b)}^{m^a}$. MacMahon's $q$-analog can also be deduced from this perspective by consideration of a \dfn{$q$-Weyl dimension formula}, that is, a principal specialization $x_i \coloneqq  q^i$ of the Weyl character formula. Indeed, whenever one has a representation-theoretic interpretation of a plane partition formula, one thus obtains a $q$-analog.

MacMahon's formula for boxed plane partitions has been generalized in many ways. Let us just mention two popular kinds of extensions in the context of lozenge tilings of the hexagon: adding \dfn{holes} to the inside of the hexagon~\cite{ciucu1998enumeration, helfhott1999enumeration, okada1998number, ciucu2001enumeration, ciucu2013dual, ciucu2017other, lai2020lozenge}; and adding \dfn{dents} to the border of the hexagon~\cite{eisenkolbl1999rhombus, ciucu2015generalization, lai2017qenumeration}. Kuo condensation is generally the most powerful tool for proving these lozenge tiling results.

\subsubsection{Shifted staircases}  \label{sec:shifted_staircase}

Set $\delta_n \coloneqq  (n,n-1,\ldots,1)$. Plane partitions in $\mathcal{SPP}^{m}(\delta_n)$, i.e., plane partitions of \dfn{shifted staircases shape}, are in bijection with symmetric plane partitions in $\mathcal{PP}^{m}(n\times n)$, as we have discussed. MacMahon~\cite{macmahon1899partitions} conjectured, and Andrews~\cite{andrews1977macmahon} proved, the following product formula enumerating symmetric plane partitions:
\[\#\mathcal{SPP}^{m}(\delta_n)= \prod_{1\leq i \leq j \leq n} \frac{m+i+j-1}{i+j-1}\]
Actually, MacMahon conjectured a $q$-analog:
\[\sum_{\substack{\pi \in \mathcal{PP}^{m}(n\times n), \\ \mathrm{Tr}(\pi)=\pi}} q^{|\pi|} =\prod_{1\leq i < j \leq n}\frac{(1-q^{2(i+j+m-1)})}{(1-q^{2(i+j-1)})} \cdot \prod_{i=1}^{n} \frac{(1-q^{2i+m-1})}{(1-q^{2i-1})}.\]
MacMahon's symmetric plane partition $q$-analog was proved by Andrews~\cite{andrews1978plane} and Macdonald~\cite{macdonald1979symmetric} essentially simultaneously and independently. Furthermore, Bender and Knuth~\cite{bender1972enumeration} conjectured a different $q$-analog for symmetric plane partitions:
\[\sum_{\substack{\pi \in \mathcal{PP}^{m}(n\times n), \\ \mathrm{Tr}(\pi)=\pi}} q^{|\pi|'} = \prod_{1\leq i \leq j \leq n}\frac{(1-q^{i+j+m-1})}{(1-q^{i+j-1})}.\]
where $|\pi|' \coloneqq  \sum_{1\leq i \leq i \leq n} \pi_{i,j}$ is a kind of ``half-size'' of $\pi$. The Bender--Knuth conjecture was proved by Gordon~\cite{gordon1983bender}, Andrews~\cite{andrews1978plane}, and Macdonald~\cite{macdonald1979symmetric}.

A representation theoretic perspective on the formula for $\mathcal{SPP}^{m}(\delta_n)$ goes as follows. Viewing the plane partitions $\pi \in \mathcal{SPP}^{m}(\delta_n)$ as Gelfand-Tsetlin patterns, one sees that they index a basis of $\oplus_{\lambda \subseteq m^n} V^{\lambda}_{\mathfrak{gl}(n)}$. Then one establishes the multiplicity free branching rule $V^{m^n}_{\mathfrak{so}(2n+1)} = \oplus_{\lambda \subseteq m^n} (\mathrm{det}^{-m/2})V^{\lambda}_{\mathfrak{gl}(n)}$ where we view $V^{m^n}_{\mathfrak{so}(2n+1)}$ as a $\mathfrak{gl}(n)$ representation under the natural inclusion $\mathfrak{gl}(n)\subseteq \mathfrak{so}(2n+1)$. The product formula then follows from the Weyl dimension formula for $V^{m^n}_{\mathfrak{so}(2n+1)}$, and the $q$-analogs follow via two slightly different principal specializations of the Weyl character formula for this $\mathfrak{so}(2n+1)$ representation. This is essentially what is done in Macdonald's proof~\cite[Examples I.5.16 and 1.5.17]{macdonald1979symmetric} (see also~\cite{okada1998applications} and~\cite{rains2015bounded} for other proofs of the $\mathfrak{so}(2n+1)$-to-$\mathfrak{gl}(n)$ branching rule). Alternatively, as Proctor observed in~\cite{proctor1984bruhat, proctor1990new}, it is possible to see directly from various known tableaux models for $\mathfrak{so}(2n+1)$ representations that $\mathcal{SPP}^{m}(\delta_n)$ indexes a basis of $V^{m^n}_{\mathfrak{so}(2n+1)}$.

For extensions of the shifted staircases formula, see the shifted double staircases discussed below.

\subsubsection{Staircases} \label{sec:staircase}

In~\cite{proctor1988odd} (see also~\cite{proctor1990new}) Proctor proved the following product formula counting plane partitions of \dfn{(unshifted) staircase shape}:
\[\#\mathcal{PP}^{m}(b,b-1,\ldots,b-a+1) =  \prod_{i=1}^{a}  \left [  \prod_{j=1}^{b-a+1} \frac{m+i+j-1}{i+j-1} \prod_{j=b-a+2}^{b-a+i} \frac{2m+i+j-1}{i+j-1} \right ].  \]
His proof uses the representation theory of the symplectic group, as we will now explain. Suppose $a+b=2N$ is even. Then Proctor~\cite{proctor1988odd} showed (essentially via King's symplectic tableaux model) that this set of plane partitions indexes a basis of~$V^{m^a}_{\mathfrak{sp}(2N)}$. The product formula follows via the Weyl dimension formula. If $a+b$ is odd, instead one has to use Proctor's ``odd symplectic group.'' It is again possible to deduce a $q$-analog via a principal specialization; however, the weights on the plane partitions that one gets are rather complicated.

The case $a=b$ corresponds to a boxed plane partition symmetry class: the plane partitions in $\mathcal{PP}^{m}(\delta_n)$ are easily seen to be in bijection with the transpose-complementary plane partition in $\mathcal{PP}^{2m}((n+1) \times (n+1))$ (see~\cite[Case 6]{stanley1986symmetries}). 

At the level of lozenge tilings, the region corresponding to the staircase plane partitions is a hexagon with ``one corner cut off.'' In~\cite{ciucu2002enumeration}, Ciucu and Krattenthaler gave an extension of Proctor's product formula to count tilings of hexagon with ``two corners cut off.''

In terms of extensions of the staircase formula, we should also mention that in~\cite{proctor1984research} (see also~\cite[Exercise~7.101]{stanley1999ec2}) Proctor himself actually established a product formula for plane partitions of a more general family of unshifted shapes which includes both the rectangles and staircases. Namely,  if we let $\lambda \coloneqq (M-d,M-2d,\ldots,M-\ell d)$ be an \dfn{arithmetic progression}, then
\[ \#\mathcal{PP}^{m}(\lambda) = \prod_{\substack{(i,j) \in \lambda, \\ \ell+c(i,j)\leq M-id} } \frac{m+\ell+c(i,j)}{\ell+c(i,j)} \cdot  \prod_{\substack{(i,j) \in \lambda, \\ \ell+c(i,j)> M-id} } \frac{(d+1)m+\ell+c(i,j)}{\ell+c(i,j)}, \]
where $c(i,j)  \coloneqq j-i$ is the \dfn{content} of the box $(i,j)\in\lambda$. The case $d=0$ corresponds to rectangles and the case $d=1$ corresponds to staircases. Proctor's proof for general $d$ is not via representation theory; rather, he observed that MacMahon's determinant in \cref{thm:macmahon_determinant} can be explicitly evaluated in this case by applying the appropriate row operations.

\subsubsection{Shifted trapezoids}

In~\cite{proctor1983trapezoid} Proctor proved the following product formula counting plane partitions of \dfn{shifted trapezoid shape}:
\[\#\mathcal{SPP}^{m}(n,n-2,\ldots,n-2(k-1)) = \prod_{i=1}^{k}\prod_{j=1}^{n-k+1}  \frac{m+i+j-1}{i+j-1}.\]
He did this again via the representation theory of the symplectic group, as we now explain. Namely, suppose $n=2N$ is even. Then Proctor~\cite{proctor1983trapezoid} showed (again using King's symplectic tableaux) that this set of plane partitions indexes a basis of~$\oplus_{\lambda \subseteq m^k} V^{\lambda}_{\mathfrak{sp}(2N)}$, and established that via the natural inclusion $\mathfrak{sp}(2N)\subseteq \mathfrak{gl}(2N+1)$ we have the multiplicity free branching rule $V^{m^k}_{\mathfrak{gl}(2N+1)}=\oplus_{\lambda \subseteq m^k} V^{\lambda}_{\mathfrak{sp}(2N)}$, from which the product formula follows via the Weyl dimension formula. The case of $n$ odd is slightly more difficult; see~\cite{proctor1983trapezoid}. Again, one can in principle obtain a $q$-analog but the weights on the plane partitions are complicated.

The case $k=\lceil n/2 \rceil$ corresponds to a boxed plane partition symmetry class: in this case the shifted trapezoid plane partitions are easily seen to be in bijection with the symmetric, self complementary plane partitions in $\mathcal{PP}^{2m}((n+1) \times (n+1))$ (see~\cite[Case 7]{stanley1986symmetries}).

In~\cite{ciucu2015lozenge}, Ciucu reproved Proctor's formula counting shifted trapezoid plane partitions via a different method: he observed that the entries in Stembridge's Pfaffian in \cref{thm:stembridge_pfaffian} radically simplify in this particular case, and the Pfaffian can be explicitly evaluated. Moreover, he in fact extended Proctor's result by allowing certain dents on the boundary of the corresponding lozenge tiling region.

\subsubsection{Shifted double staircases} 

In \cref{thm:sds} we give a product formula for plane partitions of \dfn{shifted staircase shape}, that is, the plane partitions in $\mathcal{SPP}^{m}(\lambda)$ where $\lambda \coloneqq  (n,n-1,\ldots,1) + (k,k-1,\ldots,1)$ for $0 \leq k \leq n$. Note that the case $k=0$ of this family is the shifted staircases, and the cases $k=n-1,n$ are special cases of the shifted trapezoids. Thus, this family includes the symmetric and the symmetric, self-complementary boxed plane partitions (and in some sense interpolates between them).

It is interesting to note that, unlike the case of the shifted trapezoids, the entries of Stembridge's Pfaffian in \cref{thm:stembridge_pfaffian} do not simplify for the shifted double staircases, and it thus seems hopeless to try to explicitly evaluate this Pfaffian.

Our \cref{thm:flashlight} is an extension of the shifted double staircase product formula to a more general lozenge tiling region. As mentioned, we prove this formula by extending the work of Ciucu~\cite{ciucu2020correlation}, and in particular applying the variant of Kuo condensation he developed there. One small note about deducing \cref{thm:sds} from \cref{thm:flashlight}: notice that we have the requirement $y>0$ in \cref{thm:flashlight} (see \cref{rem:y_zero}), which means that technically the case $k=n$ of \cref{thm:sds} does not follow from \cref{thm:flashlight}. However, this is okay because as mentioned the case $k=n$ of \cref{thm:sds} is a special case of Proctor's shifted trapezoid formula~\cite{proctor1983trapezoid}, and so is known.

Except for the general $d$ case of the arithmetic progression shapes, all the previously discussed shapes with product formulas have a representation theoretic interpretation. In this paper we do not give such an interpretation for the shifted double staircases. However, as mentioned in~\cref{rem:okada}, Okada's recent paper~\cite{okada2020intermediate} does give a character-theoretic extension of the shifted double staircase product formula; in particular, he obtains $q$-analogs of this product formula. Nevertheless, the approach in~\cite{okada2020intermediate} is ultimately to evaluate certain determinants and Pfaffians: a direct representation-theoretic explanation of the shifted double staircase product formula remains elusive.

\section{Proof of the main results} \label{sec:proof}

We now prove \cref{thm:flashlight} (and hence \cref{thm:sds}). We first need the following analog of Kuo condensation~\cite{kuo2004applications} for regions with a free boundary due to Ciucu~\cite{ciucu2020correlation}.

A \dfn{matching} of a graph $G$ is a collection of  vertex-disjoint edges of $G$. If a matching covers all vertices of $G$, we say that the matching is a \dfn{perfect matching} of~$G$. For a planar graph $G$ and a distinguished subset $S$ of vertices on some face, we denote by~$\M_f(G)$ the number of (not necessarily perfect) matchings $G$ in which all the vertices that are not in $S$ are matched, but those in $S$ are free to be matched or not matched (the distinguished subset of vertices $S$ will be clear from context, so we do not need to include $S$ in the notation). If $u, v, w, s \notin S$ are four vertices appearing in this cyclic order on the same face as the one containing the vertices in~$S$, we say that $S$ is `\dfn{$u, w$-separated}' if there are no mutually disjoints paths $P_1, P_2, P_3$ in $G$ so that $P_1$ connects $u$ to $w$, $P_2$ connects $v$ to some vertex in $S$, and $P_3$ connects $s$ to some other vertex of $S$. The $v,s$-separation is defined similarly.

\begin{lemma}[Corollary 1 in \cite{ciucu2020correlation}]\label{Kuolem}
Let $G$ be a planar graph with the vertices $u,v,w,s$ appearing in that cyclic order on a face $F$ of $G$. Let $S$ be a subset of the vertices of $F$ that is disjoint with $\{u,v,w,s\}$, and assume that $S$ is $u,w$- and $v,s$-separated. Then
\begin{align*}
\M_f(G)&\M_f(G\setminus\{u,v,w,s\})+\M_f(G\setminus\{u,w\})\M_f(G\setminus\{v,s\})\notag\\
&=\M_f(G\setminus\{u,s\})\M_f(G\setminus\{v,w\})+\M_f(G\setminus\{u,v\})\M_f(G\setminus\{w,s\}).
\end{align*}
\end{lemma}

We note that Ciucu actually proved a weighted version of \cref{Kuolem} in \cite{ciucu2020correlation}. However, for the purpose of our proof, we only need the above unweighted version of Ciucu's result.

\begin{figure}\centering
\includegraphics[width=4.5cm]{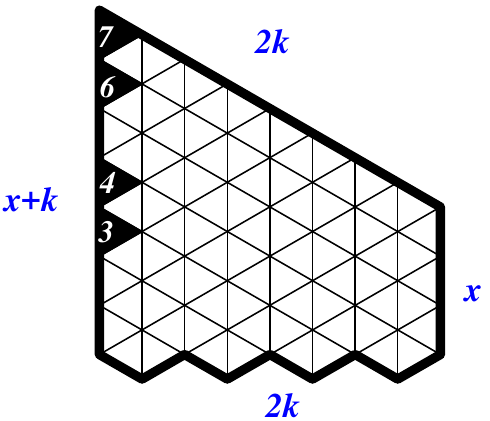}
\caption{The quartered hexagon $Q_3(3,4,6,7)$.}\label{fig:quarteredhex}
\end{figure}

We also need the following tiling enumeration of the \dfn{quartered hexagon}. A quartered hexagon $Q_{x}(s_1,s_2,\dots,s_k)$ is a trapezoidal region with sides of lengths $x+k,2k,x,2k$ as shown in \cref{fig:quarteredhex}. The left, top, and right sides of the region follow a lattice line, and the bottom side follows a zigzag lattice path. As indicated by the black triangles in that figure, we also remove $k$ right-pointing unit triangles along the left side at the positions $1\leq s_1<s_2<\cdots<s_k\leq x+k$ as they appear from the bottom to the top. The number of tilings of the quartered hexagon is always given by a simple product formula. See formula (3.2) in \cite[Theorem 3.1]{laiquarteredAztec}, or its equivalent form, formula (1.3) in \cite[Theorem 1.2]{laiquarteredhex}.\footnote{Strictly speaking, our region $Q_{x}(s_1,s_2,\dots,s_k)$ is congruent with the region $QH_{2k,x+k}(s_1,s_2,\dots,s_k)$ in \cite[Theorem 3.1]{laiquarteredAztec};  the region $Q_{x}(s_1,s_2,\dots,s_k)$ has the same tiling number as the $R_{a,b,c}$-type region (with even $b-c$) in \cite[Theorem 1.2]{laiquarteredhex}.}

\begin{lemma}[Theorem 3.1 in \cite{laiquarteredAztec}; Theorem 1.2 in \cite{laiquarteredhex}]\label{Lailem}
Assume that $x$, $k$ are non-negative integers, and that $(s_i)_{i=1}^k$ is a sequence of distinct positive integers $1\leq s_1<s_2<\cdots<s_k\leq x+k$.  The number of tilings of the quartered hexagon $Q_{x}(s_1,s_2,\dots,s_k)$ is given by
\[\M(Q_{x}(s_1,s_2,\dots,s_k))=\prod_{1\leq i < j \leq k}\frac{s_j-s_i}{j-i} \prod_{1\leq i\leq j\leq k}\frac{s_j+s_i}{j+i}.\]
\end{lemma}

We note that the number of tilings of the quartered hexagon has an interesting connection to group characters (see, e.g.,~\cite{AF, KGV}).

Each lozenge tiling of a region $R$ with a free boundary $B$ corresponds to a matching of its \dfn{(planar) dual graph} $G$ (i.e., the graph whose vertices are unit triangles in~$R$ and whose edges connect precisely two unit triangles sharing an edge), such that the vertices in $S$ correspond to unit triangles whose bases are resting on the free boundary $B$. In particular, $\M(R)=\M_f(G)$. These dual graphs are not only planar, but also \dfn{bipartite}. A bipartite graph~$G=(V_1,V_2,E)$ is a graph whose vertex set is partitioned into two vertex classes $V_1$  and $V_2$, such that two vertices in the same class are not connected by an edge. In our case, the bipartition is into \dfn{right-} and \dfn{left-pointing} triangles.

A \dfn{forced lozenge} of the region $R$ is a lozenge contained in any tiling of $R$. The removal of forced lozenges does not change the tiling number of the region.

\medskip

We are now ready to prove our main theorem.

\begin{figure}\centering
\includegraphics[width=9.5cm]{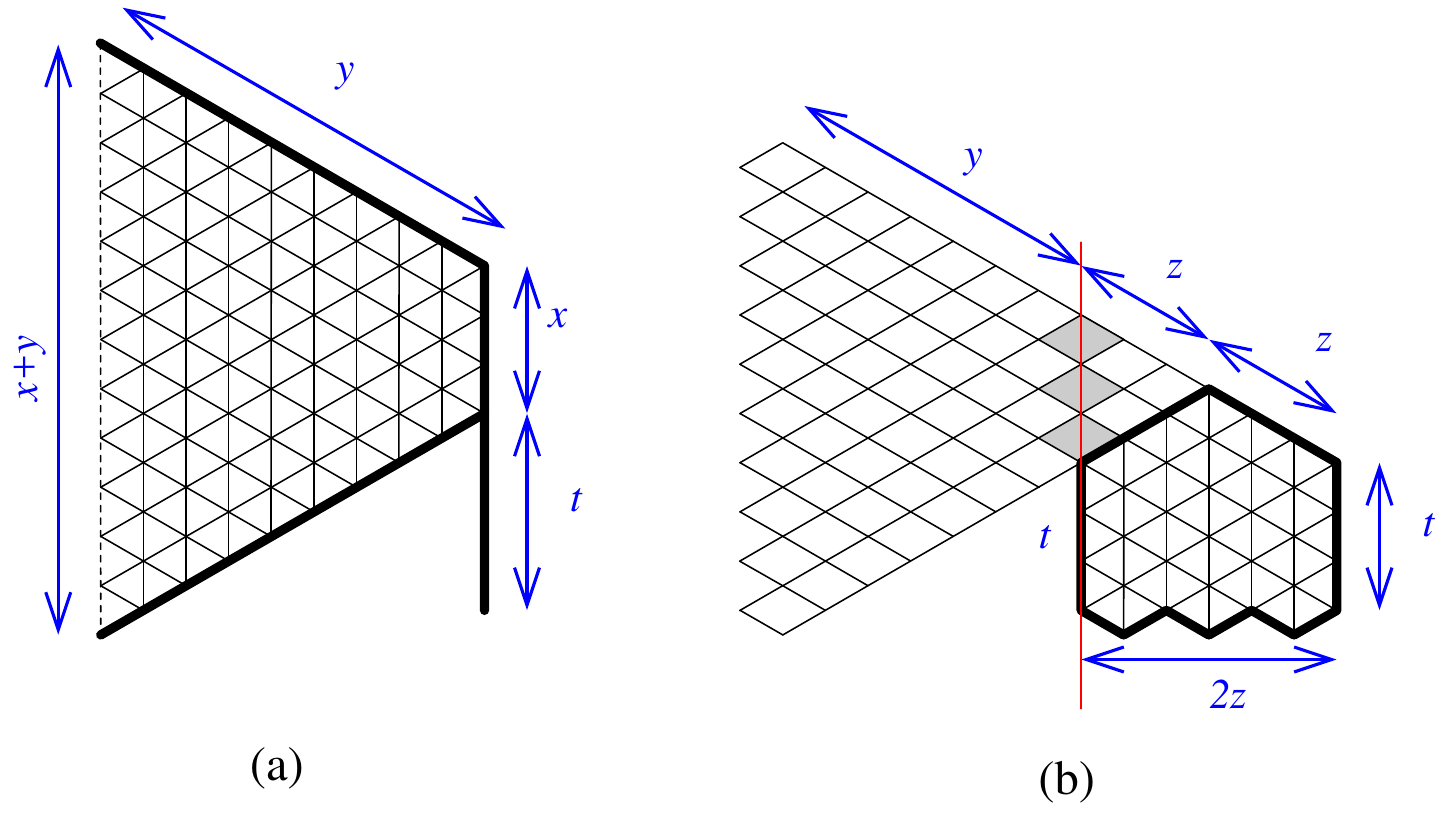}
\caption{(a) The case when $z=0$. (b) The case when $x=0$.}\label{fig:basecase1}
\end{figure}

\begin{figure}\centering
\includegraphics[width=9.5cm]{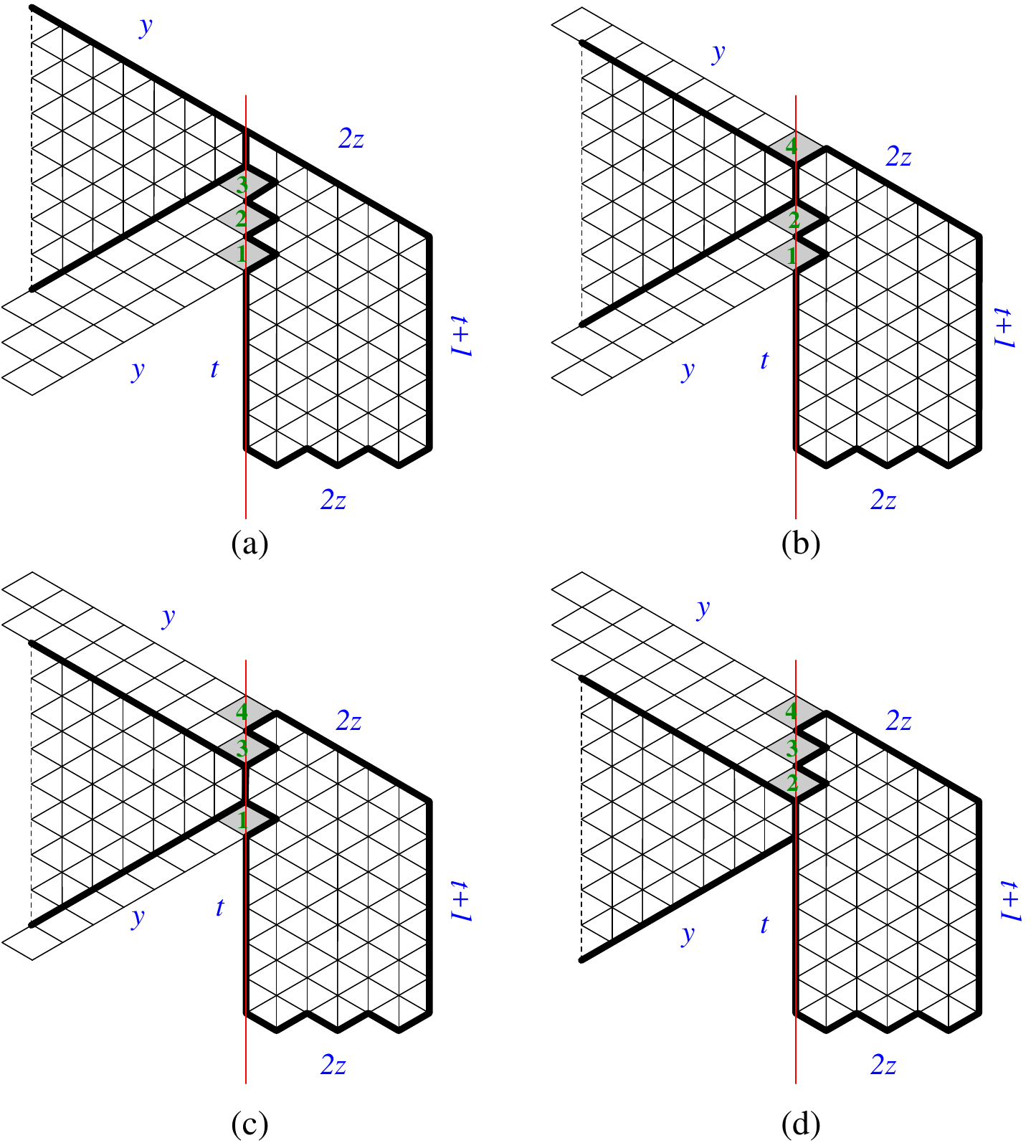}
\caption{Four ways to partition the region $F_{x,y,z,t}$ into two smaller regions in the case $x=1$.}\label{fig:basecase2}
\end{figure}

\begin{proof}[Proof of \cref{thm:flashlight}]
We prove identity~\eqref{mainformula} by induction on $x+z$. The base cases are the cases when $x=0,$ $x=1$, and $z=0$.

We first consider the case when $z=0$. In this case our region has the same tiling number as the \dfn{semi-hexagon with free boundary} (see \cref{fig:basecase1}~(a)). Thus, identity~\eqref{mainformula} follows from the well-known enumeration of symmetric plane partitions (see~\cref{sec:shifted_staircase} and~\cite{andrews1977macmahon}). 

We consider next the case when $x=0$. Let us divide the flashlight region $F_{0,y,z,t}$ along the vertical lattice line containing the $t$-side of the region (see \cref{fig:basecase1}~(b)). The right subregion has $z$ more right-pointing unit triangles than left-pointing unit triangles. It implies that, in any tiling of the right subregion, all $z$ shaded right-pointing triangles on the extended $t$-side must match outside. Equivalently, any lozenge tiling of the flashlight region must contain $z$ shaded horizontal lozenges as in the figure. The appearance of these shaded lozenges yields forced lozenges as indicated in the figure. After removing all these forced lozenges, we obtain a \dfn{halved hexagon} with the same  tiling number as $F_{0,y,z,t}$ (see the region restricted by the bold contour). The halved hexagon here is exactly the region in row 3 of Table 1 with $a=b$. The number of tilings of the halved hexagon is equal to the number of staircase-shaped plane partitions with $a=b$. These plane partitions are in bijection with transpose-complementary plane partitions (see, e.g., \cite{proctor1988odd, proctor1990new, ciucu2002enumeration}). Thus, the number of tilings of $F_{0,y,z,t}$ is equal to the number of the transpose-complementary plane partitions. However, it is not immediate to see that this agrees with our formula ~\eqref{mainformula}. We need some manipulations of the product formulas as shown below.

We define the \emph{hyperfactorial} function by: \[H(n):=\prod_{i=1}^{n} (n-i)!,\] and $H(0)=H(1)=1$. We also define the skipping version of the hyperfactorial by \[H_2(n):=\prod_{i=1}^{\lfloor n/2\rfloor}(n-2i),\] and $H_2(0)=H_2(1)=1$.

By the enumeration of transpose-complementary plane partitions and the above arguments, we have:
\begin{align}\label{identity1a}
\M(F_{0,y,z,t})&=\prod_{i=1}^{z}\frac{t+i}{i}\cdot \prod_{j=2}^{i}\frac{2t+i+j-1}{i+j-1} \notag\\
&=\prod_{i=1}^{z} \frac{2t+2i}{2i} \cdot \frac{(2t+2i-1)!/(2t+i)!}{(2i-1)!/i!} \notag\\
&= \prod_{i=1}^{z} \frac{(2t+2i)!/(2t+i)!}{(2i)!/i!}\notag\\
&=\frac{(2t+2)!(2t+4)!(2t+6)!\cdots(2t+2z)!}{2!4!6!\cdots(2z)!} \cdot \frac{1!2!3!\dots z!}{(2t+1)!(2t+2)!(2t+3)!\cdots(2t+z)!}\notag\\
&=\frac{H_2(2z+2t+2)H(z+1) H(2t+1)}{H_2(2t+2) H_2(2z+2)H(2t+z+1)}.
\end{align}

Let $P_{x,y,z,t}$ denote the expression on the right-hand side of identity~\eqref{mainformula}. We would like to verify that $\M(F_{0,y,z,t})=P_{0,y,z,t}$. It is easy to see that the first two products in the formula $P_{0,y,z,t}$ are equal to $1$, and we have
\begin{align}\label{identity1b}
P_{0,y,z,t}&=\prod_{i=1}^{z}\prod_{j=1}^{t} \frac{z+2j+i}{2j+t-1}\notag\\
&=\prod_{j=1}^{t} \frac{(2z+2j)!/(z+2j)!}{(2j+z-1)!/(2j-1)!}\notag\\
&=\frac{(2z+2)!(2z+4)!(2z+6)!\cdots(2t+2z)!}{(z+1)!(z+3)!(z+5)!\cdots (z+2t-1)!} \cdot \frac{1!3!5!\cdots (2t-1)!}{(z+2)!(z+4)!(z+6)!\cdots(z+2t)!}\notag\\
&=\frac{H_2(2t+2z+2)H_2(2t+1)H(z+1)}{H_2(2z+2)H(z+2t+1)}\notag\\
\end{align}

Note that $H(2t+1)=H_2(2t+2)H_2(2t+1)$. By \eqref{identity1a} and \eqref{identity1b}, we then have $\M(F_{0,y,z,t})=P_{0,y,z,t}$ as desired.

\bigskip

The third base case is the case when $x=1$. Arguing similarly to the case $x=0$, each lozenge tiling of the region $F_{1,y,z,t}$ must contain $z$ horizontal lozenges intersected by the extension of the $t$-side (indicated by the shaded lozenges in \cref{fig:basecase2}). Label the $z+1$ lozenges intersected by the extended $t$-side by $1,2,\dots,z+1$, from bottom to top. There are $z+1$ ways to arrange these $z$ shaded lozenges, corresponding to $z+1$ $z$-subsets $L$ of $[z+1]=\{1,2,\dots,z+1\}$. In~\cref{fig:basecase2}, we have $z=3$, and there are  four ways to arrange shaded lozenges. Let us fix an arrangement $L=\{l_1<l_2<\cdots<l_z\}$ of the $z$ shaded  lozenges and consider the region $F_{1,y,z,t}\setminus\mathcal{L}$ obtained from $F_{1,y,z,t}$ by removing these lozenges, where $\mathcal{L}$ denotes the set of  lozenges with labels in $L$. We divide this region into two smaller regions along the extended $t$-side. We can see that the right subregion is \dfn{balanced}, i.e., it has the same number of left- and right-pointing unit triangles. Thus, in any tiling of the region $F_{1,y,z,t}\setminus\mathcal{L}$, the unit triangles  in the right subregion must be matched internally. It means that any tiling of $F_{1,y,z,t}\setminus\mathcal{L}$ can be partitioned into two disjoint tilings of the left and right subregions. One could recognize  that the right subregion is congruent with a quartered hexagon $Q_{t+1}(t+l_1,t+l_2,\dots,t+l_z)$ in \cref{Lailem}. In the left subregion, there are forced lozenges along the $z$ shaded lozenges. Removing these forced lozenges, we get back  the semi-hexagon $H_s$ with free boundary.  It means that 
\[ \M(F_{1,y,z,t}\setminus\mathcal{L})=\M(H_s)\M(Q_{t+1}(t+l_1,t+l_2,\dots,t+l_z)) \]
for any $\mathcal{L}$. We note that the  shape of the semi-hexagon $H_s$ does not depend on the choice of lozenge set $\mathcal{L}$.

Summing over all $\mathcal{L}$, we get 
\[ \M(F_{1,y,z,t})=\sum_{\mathcal{L}}\M(F_{1,y,z,t}\setminus\mathcal{L})=\M(H_s)\sum_{L=\{l_1,\dots,l_z\}\subset [z+1]}\M(Q_{t+1}(t+l_1,t+l_2,\dots,t+l_z)). \]
Let us handle the sum of the tiling numbers on the right hand-side.
\begin{align*}
\sum_{\{l_1,\dots,l_z\}\subset[z+1]}\M(Q_{t+1}(t+l_1,t+l_2,&\dots,t+l_z))=\prod_{1\leq i < j \leq z}\frac{1}{j-i} \prod_{1\leq i\leq j\leq z}\frac{1}{j+i} \notag\\
&\times \left( \sum_{\{l_1,\dots,l_z\}} \prod_{1\leq i < j \leq z} (l_j-l_i) \prod_{1\leq i\leq j\leq z} (2t+l_j+l_i) \right).
\end{align*}
Compare the two products 
\[\prod_{1\leq i < j \leq z} (l_j-l_i) \prod_{1\leq i\leq j\leq z} (2t+l_j+l_i)\] and \[P=\prod_{1\leq i < j \leq z+1} (j-i) \prod_{1\leq i\leq j\leq z+1} (2t+i+j).\] 
The two products are different in the factors of $P$ that involve the term $l^*$, where $\{l^*\}=[z+1]\setminus \{l_1,\dots,l_z\}$.  By considering such factors of $P$, we get
\begin{align*}
\prod_{1\leq i < j \leq z} (l_j-l_i) \prod_{1\leq i\leq j\leq z} (2t+l_j+l_i)=P\cdot \frac{(2t+l^*)!}{(l^*-1)!(z-l^*+1)!(2t+l^*+z+1)!}.
\end{align*}
Summing over all $1\leq l^*\leq z+1$, we have
\begin{align*}
\sum_{\{l_1,\dots,l_z\}\subset[z+1]}\prod_{1\leq i < j \leq z} (l_j-l_i) &\prod_{1\leq i\leq j\leq z} (2t+l_j+l_i)\notag\\
&=P\cdot \sum_{i=1}^{z+1}\frac{(2t+i)!}{(i-1)!(z-i+1)!(2t+i+z+1)!}.
\end{align*}
Using hypergeometric series notation, we rewrite the sum on right-hand side as
\begin{align*}
\sum_{i=1}^{z+1}\frac{(2t+i)!}{(i-1)!(z-i+1)!(2t+i+z+1)!}&=\sum_{i=0}^{z}\frac{(2t+i+1)!}{(i)!(z-i)!(2t+i+z+2)!}\notag\\
&=\frac{(2t+1)!}{z!(2t+z+2)!}\sum_{i=0}^{z}\frac{(2t+2)_i(-z)_i}{(2t+z+3)_i}\frac{(-1)^i}{i!}\notag\\
&=\frac{(2t+1)!}{z!(2t+z+2)!}\cdot {}_2F_1\left[ \begin{matrix}-z & 2t+2\\2t+z+3&\end{matrix}; -1\right]
\end{align*}
By Kummer's  ${}_2F_1$ identity\footnote{Kummer's ${}_2F_1$ identity says that ${}_2F_1\left[ \begin{matrix}a & b\\c&\end{matrix}; -1\right]=\frac{\Gamma(\frac{b}{2}+1)\Gamma(b-a+1)}{\Gamma(b+1)\Gamma(\frac{b}{2}-a+1)}$ if $a-b+c=1$.}~\cite[pp. 42--43, 126]{AB}, we have 
\begin{align*}
\sum_{i=1}^{z+1}\frac{(2t+i)!}{(i-1)!(z-i+1)!(2t+i+z+1)!}&=\frac{(2t+1)!}{z!(2t+z+2)!}\frac{(t+1)!(2t+z+2)!}{(2t+2)!(t+z+1)!}\notag\\
&=\frac{t!}{2\cdot z!(t+z+1)!}.
\end{align*}
This way,  we obtain an explicit formula for the number of tilings of the flashlight region as
\begin{align*}
\M(F_{1,y,z,t})&=\M(H_s)\frac{\prod_{1\leq i < j \leq z+1} (j-i) \prod_{1\leq i\leq j\leq z+1} (2t+i+j)}{\prod_{1\leq i < j \leq z}(j-i) \prod_{1\leq i\leq j\leq z}(j+i)}\frac{t!}{2\cdot z!(t+z+1)!}\notag\\
&=\M(H_s)\frac{ \prod_{1\leq i\leq j\leq z+1} (2t+i+j)}{ \prod_{1\leq i\leq j\leq z}(j+i)}\frac{t!}{2\cdot (t+z+1)!}.
\end{align*}

By  the enumeration of symmetric plane partitions, we have
\begin{align}\label{identity2a}
\M(F_{1,y,z,t})&=\prod_{1\leq i \leq j \leq y}\frac{i+j}{i+j-1}\frac{ \prod_{1\leq i\leq j\leq z+1} (2t+i+j)}{ \prod_{1\leq i\leq j\leq z}(j+i)}\frac{t!}{2\cdot (t+z+1)!}\notag\\
&=2^y\frac{ \prod_{1\leq j\leq z+1} (2t+2j)!/(2t+j)!}{ \prod_{1\leq j\leq z}(2j)!/j!}\frac{t!}{2\cdot (t+z+1)!}\notag\\
&=\frac{2^{y-1}t!}{(t+z+1)!}\frac{H_2(2t+2z+4)H(2t+1)H(z+1)}{H_2(2t+2)H(2t+z+2)H_2(2z+2)}.
\end{align}

We would like to verify that $\M(F_{1,y,z,t})=P_{1,y,z,t}$. We can simplify the expression $P_{1,y,z,t}$ on the right-hand side of  \eqref{mainformula} as:
\begin{align}\label{identity2b}
P_{1,y,z,t}&=\prod_{1\leq j \leq y+z}\frac{2j}{j}\prod_{1\leq j \leq z} \frac{2j+1}{j+1} \prod_{i=1}^{t}\prod_{j=1}^{z}\frac{2i+j+z+1}{2i+j}\notag\\
&=2^{y+z}\frac{(2z+1)!!}{(z+1)!}\prod_{i=1}^{t}\frac{(2i+2z+1)!/(2i+z+1)!}{(2i+z)!/(2i)!}\notag\\
&=\frac{2^{y+z}(2z+1)!!}{(z+1)!}\frac{H_2(2z+2t+3)H_2(2t+2)H(z+1)}{H_2(2z+3)H(2t+z+2)}.
\end{align}

We note that \[H_2(2x+2)=(2x)!! H_2(2x+1)=2^{x}x! H_2(2x+1)\] and \[H_2(2x+1)=(2x-1)!! H_2(2x).\] This and the above two identities \eqref{identity2a} and  \eqref{identity2b} imply that $\M(F_{1,y,z,t})=P_{1,x,y,t}$ as desired. Then identity~\eqref{mainformula} follows for the case $x=1$.

\medskip

For the induction step we assume that $x\geq 2$, $z\geq 1$ and that~\eqref{mainformula} holds for any flashlight region whose sum of $x$- and $z$-parameter is strictly less than $x+z$. We will use Ciucu's~\cref{Kuolem} to construct certain recurrences for the tiling numbers of the flashlight regions, and~\eqref{mainformula} follows from the induction principle.

\begin{figure}\centering
\includegraphics[width=12cm]{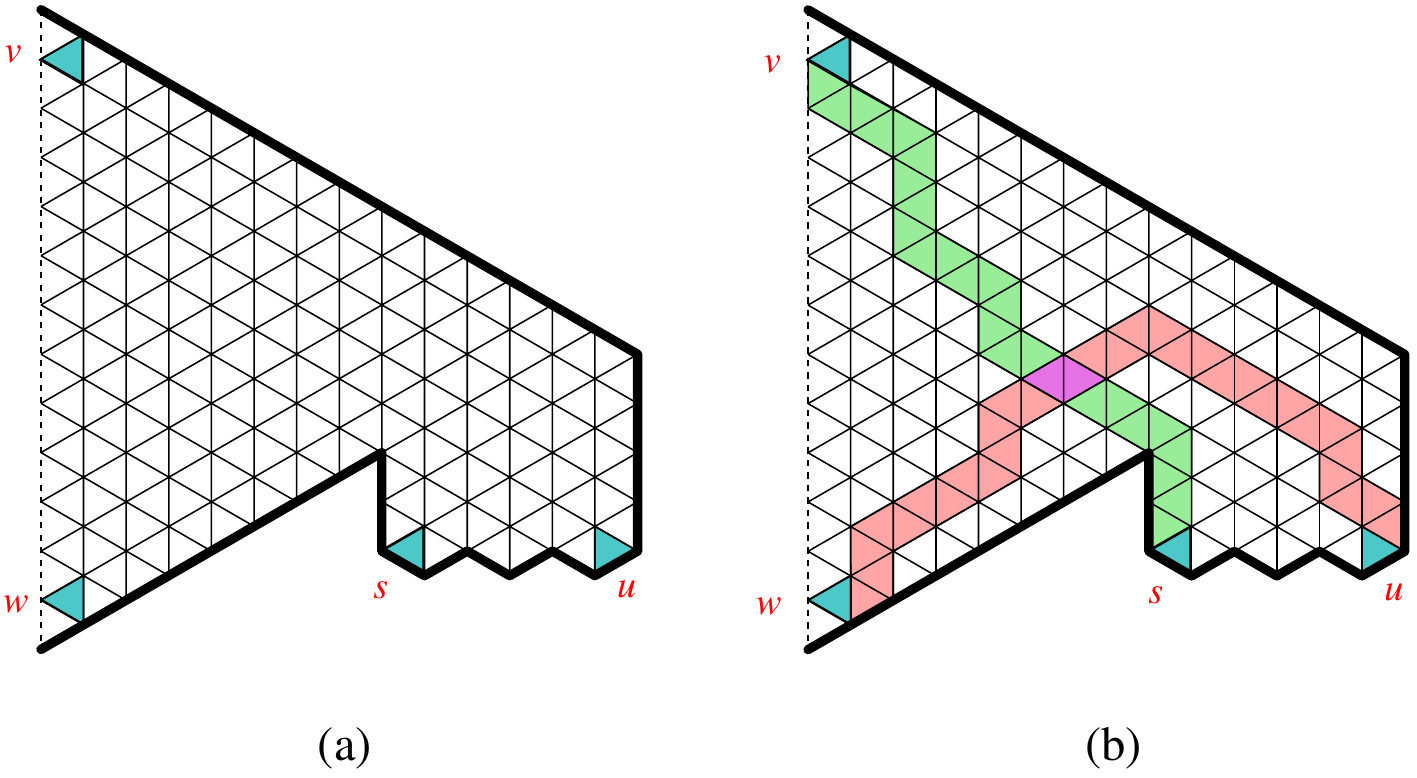}
\caption{(a) How to apply the condensation in \cref{Kuolem} to the flashlight region. (b) Separating condition in Ciucu's~\cref{Kuolem}}\label{fig:Kuo1}.
\end{figure}

\begin{figure}\centering
\includegraphics[width=9.75cm]{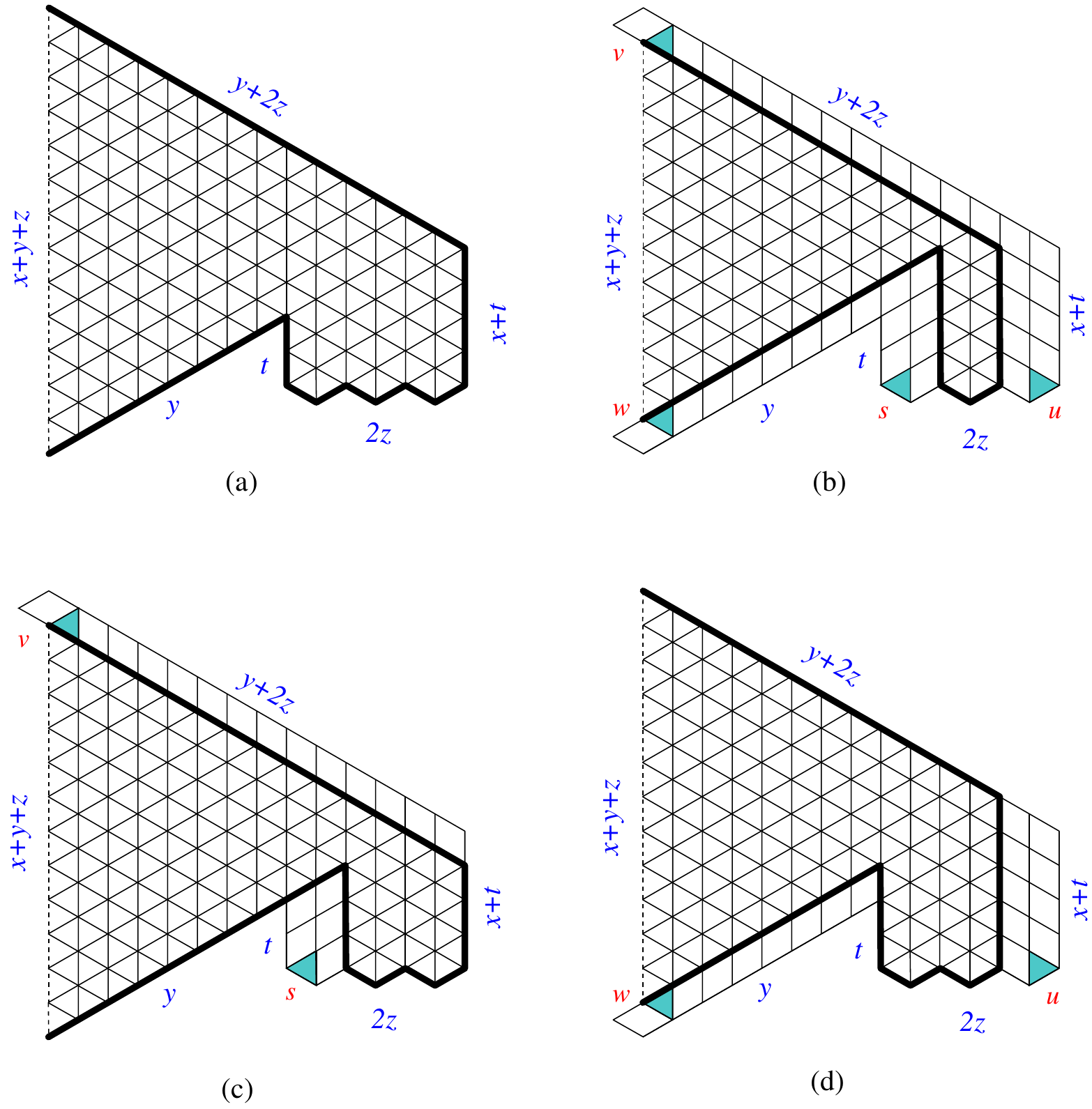}
\caption{The four regions  corresponding to the four graphs on the left-hand side of identity~\eqref{kuoeqb}.}\label{fig:Kuo2}
\end{figure}

\begin{figure}\centering
\includegraphics[width=9.75cm]{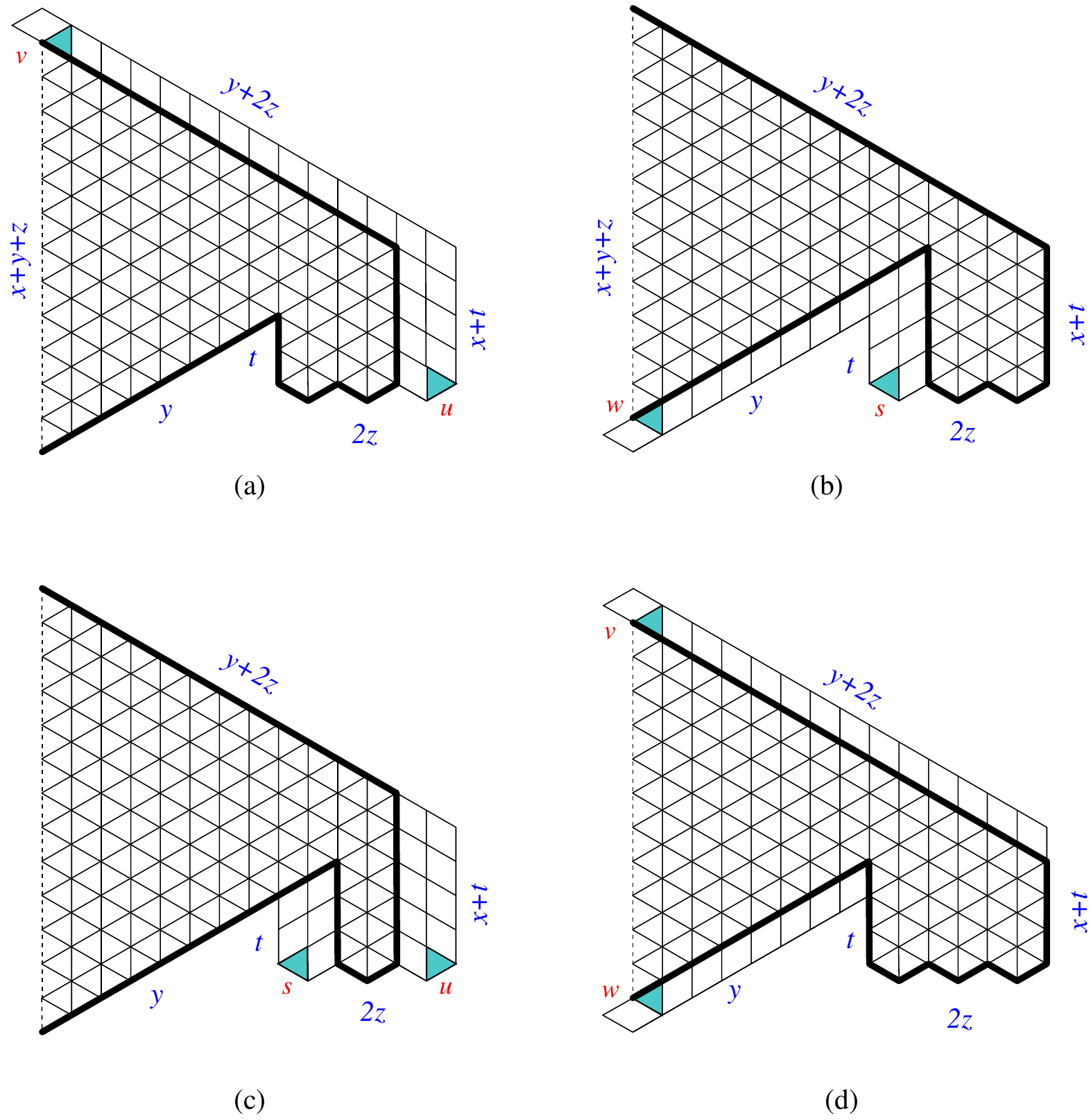}
\caption{The four regions  corresponding to the four graphs on the right-hand side of identity~\eqref{kuoeqb}.}\label{fig:Kuo3}
\end{figure}

We first suppose that $z\geq 2$. Apply Ciucu's \cref{Kuolem} to the dual graph $G$ of the flashlight region $F_{x,y,z,t}$ with the four vertices $u,v,w,s$ chosen as in \cref{fig:Kuo1}(a). Strictly speaking, the vertices $u,v,w,s$ of $G$ are indicated by the corresponding unit triangles in the region $F_{x,y,z,t}$ (see the shaded unit triangles of the same label). In particular, the $u$-triangle is the right-pointing unit triangle in the lower-right corner of the region, the $v$-, $w$-, and $s$-triangles  are the shaded ones as we go in counter-clockwise direction around the region from the $u$-triangle. The separating conditions in Ciucu's Lemma is now equivalent to the fact that $s$ and $v$ belong to different parts of the region $F_{x,y,z,t}$ separated by any path of triangles going from $u$ to $w$, and $u$ and $w$ belong to different parts of the region  $F_{x,y,z,t}$ separated by any path of triangles going from $s$ to $v$ (see \cref{fig:Kuo1}(b)). We now have the recurrence:
\begin{align}\label{kuoeqb}
\M_f(G)\M_f&(G\setminus\{u,v,w,s\})+\M_f(G\setminus\{v,s\})\M_f(G\setminus\{u,w\})\notag\\
&=\M_f(G\setminus\{u,v\})\M_f(G\setminus\{w,s\})+\M_f(G\setminus\{u,s\})\M_f(G\setminus\{v,w\}).
\end{align}

We would like to convert this recurrence to a recurrence for tiling numbers of the flashlight regions. We will write each of the eight matching numbers in terms of tiling number of some flashlight region. The task for the first matching number is trivial, as we have by definition:
\begin{equation}\label{convert1}
\M_f(G)=\M(F_{x,y,z,t}).
\end{equation}
Let us consider the region corresponding to the graph in  the second matching number, i.e., the graph $G\setminus\{u,v,w,s\}$. The removal of  the $u$-, $v$-, $w$-, and $s$-triangles yields several forced lozenges as indicated as in \cref{fig:Kuo2}(b). We note that the removal of the $v$- and $w$-triangles forces the  topmost and bottommost  right-pointing triangles along the free boundary to match outside. Removing these forced lozenges, we get a new flashlight region (indicated by the region restricted by the bold contour), namely $F_{x-2,y+2,z-2,t+2}$. As removal of forced lozenges does not change the tiling number, we have:
\begin{equation}\label{convert2}
\M_f(G\setminus\{u,v,w,s\})=\M(F_{x-2,y+2,z-2,t+2}).
\end{equation}
Considering forced lozenges as shown in \cref{fig:Kuo2}(c), (d) and \cref{fig:Kuo3}(a)--(d), we get respectively six more equations:
\begin{equation}\label{convert3}
\M_f(G\setminus\{v,s\})=\M(F_{x-2,y+2,z-1,t+1}),
\end{equation}
\begin{equation}\label{convert4}
\M_f(G\setminus\{u,w\})=\M(F_{x,y,z-1,t+1}),
\end{equation}
\begin{equation}\label{convert5}
\M_f(G\setminus\{u,v\})=\M(F_{x,y,z-1,t}),
\end{equation}
\begin{equation}\label{convert6}
\M_f(G\setminus\{w,s\})=\M(F_{x-2,y+2,z-1,t+2}),
\end{equation}
\begin{equation}\label{convert7}
\M_f(G\setminus\{u,s\})=\M(F_{x,y+2,z-2,t+1}),
\end{equation}
\begin{equation}\label{convert8}
\M_f(G\setminus\{v,w\})=\M(F_{x-2,y,z,t+1}).
\end{equation}
Plugging equations~\eqref{convert1}--\eqref{convert8} in recurrence~\eqref{kuoeqb}, we get a recurrence for tiling numbers of flashlight regions:
\begin{align}\label{kuoeqc}
\M(&F_{x,y,z,t})\M(F_{x-2,y+2,z-2,t+2})+\M(F_{x-2,y+2,z-1,t+1})\M(F_{x,y,z-1,t+1})\notag\\
&=\M(F_{x,y,z-1,t})\M(F_{x-2,y+2,z-1,t+2})+\M(F_{x,y+2,z-2,t+1})\M(F_{x-2,y,z,t+1}).
\end{align}
One readily sees that the sums of $x$- and $z$-parameters of the last seven regions in the above recurrence are all strictly less than $x+z$. By induction hypothesis, we have explicit formulas for the last seven tiling numbers.  Then we obtain the formula for the number of tilings of $F_{x,y,z,t}$, which is the same as the expression on the right-hand side of identity~\eqref{mainformula}, say after performing a straightforward simplification. Recall that $P_{x,y,z,t}$ denotes the expression on the right-hand side of identity~\eqref{mainformula}. We remark that not only do these $P_{x,y,z,t}$ satisfy the recurrence~\eqref{kuoeqc}, but in fact it could be seen that
\begin{equation}\label{Pidentity1} P_{x,y,z,t}P_{x-2,y+2,z-2,t+2} = P_{x,y+2,z-2,t+1}P_{x-2,y,z,t+1}\end{equation}
and
\begin{equation}\label{Pidentity2} P_{x-2,y+2,z-1,t+1}P_{x,y,z-1,t+1} = P_{x,y,z-1,t}P_{x-2,y+2,z-1,t+2}.\end{equation}
Again, verifying these identities is not hard but not immediate either. We will present briefly the verification below.

Identity~\eqref{Pidentity2} is easier to see by showing that
\begin{align}
\frac{P_{x-2,y+2,z-1,t+1}}{P_{x-2,y+2,z-1,t+2}}\cdot \frac{P_{x,y,z-1,t+1} }{P_{x,y,z-1,t}}=1.
\end{align}
We note that in each of the fractions on the left-hand side, the two $P$-terms are the same, except for the $t$-parameters. Cancelling out all factors independent from the $t$-parameter, we have
\begin{align}
\frac{P_{x-2,y+2,z-1,t+1}}{P_{x-2,y+2,z-1,t+2}}\cdot \frac{P_{x,y,z-1,t+1} }{P_{x,y,z-1,t}}&=\frac{\prod_{i=1}^{t+1}\prod_{j=1}^{z-1}\frac{x+z+2i+j-3}{x+2i+j-3}}{\prod_{i=1}^{t+2}\prod_{j=1}^{z-1}\frac{x+z+2i+j-3}{x+2i+j-3}} \cdot\frac{\prod_{i=1}^{t+1}\prod_{j=1}^{z-1}\frac{x+z+2i+j-1}{x+2i+j-1}}{\prod_{i=1}^{t}\prod_{j=1}^{z-1}\frac{x+z+2i+j-1}{x+2i+j-1}}\notag\\
&=\frac{1}{\prod_{j=1}^{z-1}\frac{x+z+2t+j+1}{x+2t+j+1}}\cdot \frac{\prod_{j=1}^{z-1}\frac{x+z+2t+j+1}{x+2t+j+1}}{1}\notag\\
&=1.
\end{align}

Next, we prove identity~\eqref{Pidentity1} by verifying that 
\begin{equation}\frac{ P_{x,y,z,t}P_{x-2,y+2,z-2,t+2}}{ P_{x,y+2,z-2,t+1}P_{x-2,y,z,t+1}}=1.\end{equation}
Cancelling common terms in the numerator and denominator on the left-hand side, we have
\begin{align}\frac{ P_{x,y,z,t}P_{x-2,y+2,z-2,t+2}}{ P_{x,y+2,z-2,t+1}P_{x-2,y,z,t+1}}&=\frac{\prod_{1\leq i \leq j \leq z}x+i+j \prod_{1\leq i \leq j \leq z-2}x+i+j-2}{\prod_{1\leq i \leq j \leq z-2}x+i+j \prod_{1\leq i \leq j \leq z}x+i+j-2}\notag\\
&\quad\times\frac{\prod_{i=1}^{t}\prod_{j=1}^{z}\frac{x+z+2i+j}{x+2i+j-1} \prod_{i=1}^{t+2}\prod_{j=1}^{z-2}\frac{x+z+2i+j-4}{x+2i+j-3}}{\prod_{i=1}^{t+1}\prod_{j=1}^{z-2}\frac{x+z+2i+j-2}{x+2i+j-1}\prod_{i=1}^{t+1}\prod_{j=1}^{z}\frac{x+z+2i+j-2}{x+2i+j-1}}\notag\\
&=\frac{\prod_{1\leq i \leq z-1}x+z+i-1 \prod_{1\leq i \leq z}x+z+i}{\prod_{1\leq i \leq z-1}x+z+i-3 \prod_{1\leq i \leq z}x+z+i-2}\notag\\
&\quad\times \prod_{j=1}{z}\frac{x+j-1}{x+z+j}\prod_{i=1}^{z-2}\frac{x+z+j-2}{x+j-1} \notag\\
&=\frac{(x+2z-3)(x+2z-2)}{(x+z-3)(x+z-1)} \frac{(x+2z)(x+2z-1)}{(x+z-1)(x+z-2)} \notag\\
&\quad\times\frac{(x+z-3)(x+z-1)(x+z-1)(x+z)}{(x+2z)(x+2z-1)(x+2z-2)(x+2z-3)}\notag\\
&=1.
\end{align}

We have just proved identity~\eqref{mainformula} for the case $z\geq 2$.

\medskip

\begin{figure}\centering
\includegraphics[width=8.75cm]{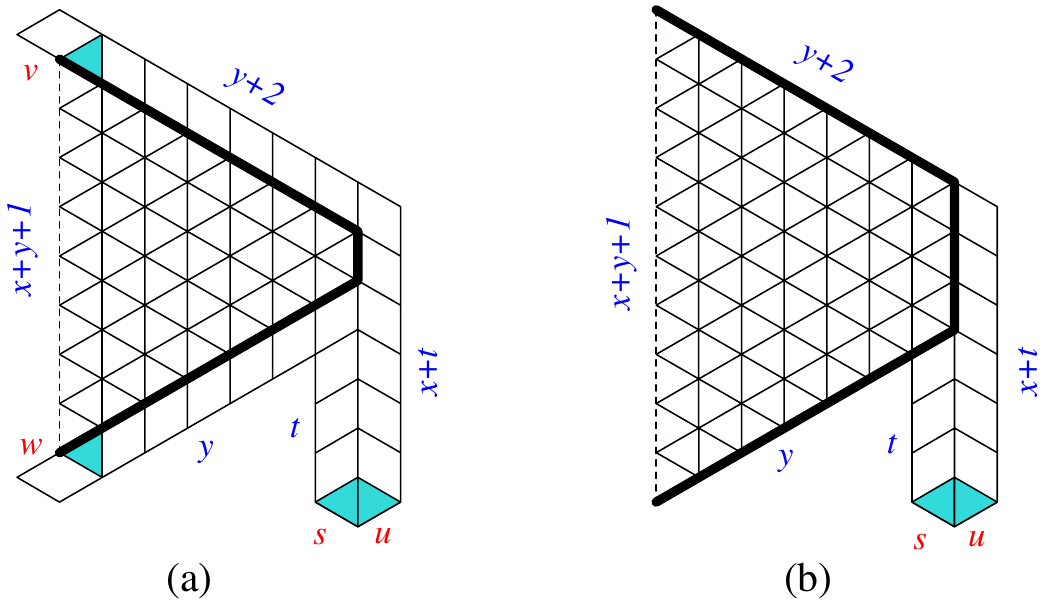}
\caption{Two special regions in the case when $x\geq 2$ and $z=1$.}\label{fig:Kuo4}
\end{figure}

We note that recurrence~\eqref{kuoeqc} does not work for the case $z=1$. In particular, the regions $F_{x-2,y+2,z-2,t+2}$ and $F_{x,y+2,z-2,t+1}$ are not defined when $z=1$. To complete the proof, we need to obtain a new recurrence in this case. We still apply Ciucu's \cref{Kuolem} to the region $F_{x,y,1,t}$ as in \cref{fig:Kuo1}. It is easy to see that among the above eight equations~\eqref{convert1}--\eqref{convert8}, only two equations do not hold when $z=1$. They  are equations~\eqref{convert2} and~\eqref{convert7}. We now need to consider the effect of forced lozenges in the regions corresponding to the graphs $G\setminus\{u,v,w,s\}$ and $G\setminus\{u,s\}$ when $z=1$, as shown in \cref{fig:Kuo4} (a) and (b), respectively. In particular, we get two new equations:
\begin{equation}\label{convert2b}
\M_f(G\setminus\{u,v,w,s\})=\M(F_{x-2,y+1,0,t+2})
\end{equation}
and
\begin{equation}\label{convert7b}
\M_f(G\setminus\{u,s\})=\M(F_{x,y+1,0,t+1}).
\end{equation}
Using~\eqref{convert2b} and~\eqref{convert7b} instead of~\eqref{convert2} and~\eqref{convert7}, we get a new recurrence for the case $z=1$:
\begin{align}\label{kuoeqd}
\M(&F_{x,y,1,t})\M(F_{x-2,y+1,0,t+2})+\M(F_{x-2,y+2,0,t+1})\M(F_{x,y,0,t+1})\notag\\
&=\M(F_{x,y,0,t})\M(F_{x-2,y+2,0,t+2})+\M(F_{x,y+1,0,t+1})\M(F_{x-2,y,1,t+1}).
\end{align}
(We note that, as mentioned in the case $z=0$, the tiling number of the region $F_{x,y,0,t}$ is actually independent from $t$.)
Again, the sums of the $x$- and $z$-parameters of the last seven regions in the recurrence~\eqref{kuoeqd} are all smaller than $x+z$, and~\eqref{mainformula} also follows from the induction hypothesis in this case.  As before, letting $P_{x,y,z,t}$ denote the expression on the right-hand side of identity~\eqref{mainformula}, we remark that not only do these $P_{x,y,z,t}$ satisfy the recurrence~\eqref{kuoeqd}, but in fact
\[P_{x,y,1,t}P_{x-2,y+1,0,t+2} = P_{x,y+1,0,t+1}P_{x-2,y,1,t+1},\]
and 
\[P_{x-2,y+2,0,t+1}P_{x,y,0,t+1} = P_{x,y,0,t}P_{x-2,y+2,0,t+2},\]
which can be seen by similar computations to what we showed in the case $z\geq 2$. This completes our proof.
\end{proof}

\begin{remark} \label{rem:y_zero}
Although we stated \cref{thm:flashlight} for $y>0$, we believe it holds verbatim in the case $y=0$ as well. The problem is that we cannot use \cref{Kuolem} in this case, because the relevant graph is no longer $u, w$-separated. The number of tilings $\M(F_{x,0,z,t})$ is naturally expressed as a sum of $\M(Q_k(s_1,\ldots,s_k))$, but this sum seems hard to evaluate.
\end{remark}

\bibliography{sds}{}
\bibliographystyle{plain}

\end{document}